\pdfoutput=1
\RequirePackage{ifpdf}
\ifpdf % We~are running pdfTeX in pdf mode
\documentclass[pdftex]{sigma}%,draft
\else
\documentclass{sigma}
\fi

\usepackage{cleveref}

\usepackage{bbm}
\newcommand{\Eins}{\mathbbm{1}}

\usepackage{enumitem}
\usepackage{tikz-cd} % for commutative diagrams
\usepackage{stmaryrd}
\usepackage{mathrsfs}

\newcommand{\id}{\mathrm{id}}
\newcommand{\betrag}[1]{\left|#1\right|}
\newcommand{\Alg}{\mathfrak{alg}}
\newcommand{\algQ}{\mathfrak{algQ}}
\newcommand{\Vect}{\mathfrak{vec}}
\newcommand{\Hilb}{\mathfrak{hilb}}
\newcommand{\Set}{\mathfrak{set}}
\newcommand{\Prob}{\mathfrak{prob}}
\newcommand{\FinVect}{\mathfrak{finvec}}
\newcommand{\FinHilb}{\mathfrak{finhilb}}

\newcommand{\Mor}{\operatorname{\mathfrak{mor}}}
\newcommand{\Obj}{\operatorname{\mathfrak{obj}}}

\newtheorem{Theorem}{Theorem}[section]
\newtheorem{Corollary}[Theorem]{Corollary}
\newtheorem{Lemma}[Theorem]{Lemma}
\newtheorem{Proposition}[Theorem]{Proposition}
\newtheorem{Observation}[Theorem]{Observation}
\newtheorem*{Theorem*}{Theorem}

{ \theoremstyle{definition}
\newtheorem{Definition}[Theorem]{Definition}

\newtheorem{Example}[Theorem]{Example}
\newtheorem{Remark}[Theorem]{Remark}}

\crefname{Definition}{Definition}{Definitions}
\crefname{Example}{Example}{Examples}
\crefname{Note}{Note}{Notes}
\crefname{Remark}{Remark}{Remarks}

\crefname{subsection}{subsection}{subsections}

\numberwithin{equation}{section}

\begin{document}
\allowdisplaybreaks

\newcommand{\arXivNumber}{1612.05139}

\renewcommand{\thefootnote}{}

\renewcommand{\PaperNumber}{075}

\FirstPageHeading

\ShortArticleName{Categorial Independence and L\'evy Processes}

\ArticleName{Categorial Independence and L\'evy Processes\footnote{This paper is a~contribution to the Special Issue on Non-Commutative Algebra, Probability and Analysis in Action. The~full collection is available at \href{https://www.emis.de/journals/SIGMA/non-commutative-probability.html}{https://www.emis.de/journals/SIGMA/non-commutative-probability.html}}}

\Author{Malte GERHOLD~$^{\rm ab}$, Stephanie LACHS~$^{\rm a}$ and Michael SCH\"URMANN~$^{\rm a}$}

\AuthorNameForHeading{M.~Gerhold, S.~Lachs and M.~Sch\"urmann}

\Address{$^{\rm a)}$~Institute of Mathematics and Computer Science, University of Greifswald, Germany}
\EmailD{\href{mailto:gerholds@uni-greifswald.de}{gerholds@uni-greifswald.de}, \href{mailto:schurman@uni-greifswald.de}{schurman@uni-greifswald.de}}

\Address{$^{\rm b)}$~Department of Mathematical Sciences, NTNU Trondheim, Norway}
\EmailD{\href{mailto:malte.gerhold@ntnu.no}{malte.gerhold@ntnu.no}}

\ArticleDates{Received March 28, 2022, in final form September 30, 2022; Published online October 10, 2022}

\Abstract{We generalize Franz' independence in tensor categories with inclusions from two morphisms (which represent generalized random variables) to arbitrary ordered families of morphisms. We will see that this only works consistently if the unit object is an initial object, in which case the inclusions can be defined starting from the tensor category alone. The obtained independence for morphisms is called categorial independence. We define categorial L\'evy processes on every tensor category with initial unit object and present a~construction generalizing the reconstruction of a L\'evy process from its convolution semigroup via the Daniell--Kolmogorov theorem. Finally, we discuss examples showing that many known independences from algebra as well as from (noncommutative) probability are special cases of categorial independence.}


\Classification{18D10; 60G20; 81R50}

\renewcommand{\thefootnote}{\arabic{footnote}}
\setcounter{footnote}{0}

\section{Introduction}\label{sec:introduction}

Suppose $(\mu_t)_{t\in\mathbb R_+}$ is a convolution semigroup of probability measures on the real line. Let us sketch, how to construct a L\'evy process $X_t\colon\Omega\to\mathbb R$ with marginal distributions $\mathbb{P}_{X_t}=\mu_t$. First, for all finite subsets $J=\{t_1<t_2<\cdots <t_n\}\subset\mathbb R_+$ define probability measures $\mu_J:=\mu_{t_1}\otimes\mu_{t_2-t_1}\otimes\cdots\otimes\mu_{t_n-t_{n-1}}$ on $\mathbb R^J$ . Then show that these are coherent in the sense that
\begin{gather*}
 \mu_I=\mu_J\circ \big(p^{J}_I\big)^{-1}
\end{gather*}
for all $I\subset J$ and $p^{J}_I\colon \mathbb R^J\to\mathbb R^I$ the canonical projection. In this situation the probability spaces $\big(\mathbb R^J,\mathcal B\big(\mathbb R^J\big),\mu_J\big)$ with the projections $\big(p^{J}_I\big)$ form a \emph{projective system}. Now, the Daniell--Kolmogorov theorem guarantees the existence of a \emph{projective limit}, which is a probability space $(\Omega,\mathcal F,\mathbb{P})$ with projections $p_J\colon \Omega\to \mathbb R^{J}$ such that $\mu_J=\mathbb{P}\circ(p_J)^{-1}$. The random variable ${X_t:=p_{\{t\}}}$ has distribution $\mu_t$ and the $X_t$ have independent and stationary increments ${X_s-X_t\sim\mu_{s-t}}$.

A similar construction allows one to associate quantum L\'evy processes with convolution semigroups of states on $*$-bialgebras. The formulation of quantum probability is dual to that of classical probability, so inductive limits appear instead of projective limits. Due to the fact that there are different independences in quantum probability on the one hand and the interactions between quantum probability and operator algebras on the other hand, there are many different theorems of the same kind (construction of L\'evy processes for other notions of independence, see \cite{bGSc05, Sch93}) or similar kind (the construction of product systems and dilations from subproduct systems, cf.\ \cite{BhMu10, BhSk00}). Also, following Voiculescu's invention of bifreeness~\cite{Voi14}, many examples of multivariate independences have been exhibited \cite{G17p,GHU21p,GHS20,GuSk17,GuSk19,Liu18p,Liu19} and some general theory has been developed~\cite{G21p, MaSc17}, and all those independences have associated classes of L\'evy processes that are very interesting to study. The main aim of these notes is to give a unified approach to these different situations. To this end, we work with the language of tensor categories and introduce \emph{comonoidal systems}, which serve as a suitable generalization of convolution semigroups.

With a similar goal in mind, Uwe Franz defined independence on \emph{tensor categories with inclusions} \cite{Fra06}.\footnote{Following Franz, we use the terms \emph{tensor category}, \emph{tensor functor}, and \emph{cotensor functor} synonymously with \emph{monoidal category}, \emph{lax monoidal functor}, and \emph{colax monoidal functor}, respectively, see \Cref{sec:indep-tens-categ}; while the latter are more common in abstract category theory, the former are frequently used in applications such as noncommutative probability or quantum group theory.} Whereas his definition works nicely for all questions concerning independence of two random variables, the approach runs into trouble when there are more random variables involved, which inevitably happens for example when one wants to study general properties of L\'evy processes with respect to different notions of independence. We solve this problem by considering tensor categories with initial unit object instead of tensor categories with inclusions. We show that inclusions in the sense of Franz can always be defined in such a category and determine the necessary and sufficient conditions for given inclusions to be of the desired form.\footnote{Of course, tensor categories with terminal unit object and the corresponding projections can be also be considered and all results apply dually (i.e., with arrows reversed); those are also known to category theorists under the name \emph{semicartesian categories} and play an important role in the synthetic approach to probability using Markov categories, see \cite{Fri20} and references therein.} In this stronger setting, Franz' notion of independence can be generalized to arbitrary ordered families of random variables and will be called \emph{categorial independence}. Based on categorial independence, we will define categorial L\'evy processes and present a general construction of L\'evy processes from comonoidal systems reminescent of the reconstruction of a L\'evy process from its convolution semigroup via the Daniell--Kolmogorov theorem.

We will see that categorial independence is a very general concept of independence which encompasses most independences encountered in different areas of mathematics, such as linear independence, algebraic independence, orthogonality or stochastic independence.

A widely known approach to general independences is via matroids. This is very different from the approach treated in this article and, indeed, it is quite obvious that the systems of orthogonal sets in a Hilbert space, or of independent sets of random variables in a probability space do not form a matroid, while they turn out to be special cases of categorial independence.

Another very general concept of independence, somehow closer to categorial independence, is due to Marczewski \cite{Mar58} and is explored to some detail in Gr{\"a}tzer's book on universal algebra \cite{Grae79}. This is a concept of independence meant to generalize and unify independences mainly arising in algebra, in particular, Marczewski clarifies that stochastic independence cannot be obtained this way \cite[Section 9]{Mar58}. Therefore, categorial independence is not a special case of Marczewski's independence. Conversely, our list of examples includes Marczewski's examples from algebra. However, we leave as an open problem the question whether Marczewski's independence and other related notions from universal algebra can be covered by categorial independence in full generality. In particular, we would be very interested to know whether logical independence is a special case of categorial independence.

Independently of this work, Simpson \cite{Sim18} recently developed a categorical approach to independence that seems closely related to ours; in many cases the tensor product in our approach is the same as the \emph{independent product} in his.

The paper is organized as follows. In \Cref{sec:tensor-categories} we recall basic facts concerning inductive limits and tensor categories. In \Cref{sec:indep-tens-categ} we recall Franz' original definition of categorial independence in tensor categories with inclusions and generalize it to more than two morphisms in case the inclusions are compatible; also, we prove that inclusions are compatible if and only if they are the canonical inclusions in a tensor category with initial unit object. In \Cref{sec:comonoidal-systems} we introduce a categorial counterpart of convolution semigroups, which we call \emph{comonoidal systems} and prove that under mild assumptions comonoidal systems give rise to categorial L\'evy processes by inductive limit constructions. Finally, in \Cref{sec:examples-1} we discuss examples from algebra and (noncommutative) probability.

\section{Basic notions of category theory}\label{sec:tensor-categories}

The main point of this section is to fix notations and recall basic facts about inductive limits and tensor categories. We also give a list of those categories which appear as examples for the following sections.

We will freely use the language of categories, functors, and natural transformations, see for example the book of Ad{\'a}mek, Herrlich and Strecker \cite{AHS04} for an exquisite treatment of the matter. For a category $\mathcal{C}$, we write $\Obj \mathcal{C}$ for the class of objects of $\mathcal C$. Given two objects $A,B\in\Obj\mathcal{C}$, we denote by $\Mor(A,B)$ the set of all morphisms $f\colon A\to B$.

Equalities between morphisms will frequently be expressed in terms of commutative diagrams. A \emph{diagram} is a directed graph with object-labeled vertices and morphism-labeled edges. We say that a diagram \emph{commutes} if the composition of morphisms along any two directed paths with the same source and the same target vertex yield the same result. We will usually not explicitly write the inverse of an isomorphism with an extra edge, but it shall be included when we say that the diagram commutes. If the morphism labelling an edge is a component $\alpha_A$ of a natural transformation $\alpha$, and it is clear from the context (i.e., the source and the target object) which component we mean, we will drop the index to increase readibility.

\subsection{Inductive limits}\label{sec:inductive-limits}

In category theory there are the general concepts of limits and colimits. Since in our applications only inductive limits play a role, we restrict to this special case. The general case can for example be found in~\cite{AHS04}.

A preordered set $I$ is called \emph{directed} if any two elements of $I$ possess a common upper bound, that is if for all $\alpha,\beta\in I$ there exists $\gamma\in I$ with $\gamma\geq\alpha,\beta$.

\begin{Definition} Let $\mathcal C$ be a category. An \emph{inductive system} consists of
 \begin{itemize}\itemsep=0pt
 \item a family of objects $(A_\alpha)_{\alpha\in I}$ indexed by a directed set $I$,
 \item a family of morphisms $(f^{\alpha}_{\beta}\colon A_\alpha\to A_\beta)_{\alpha\leq\beta}$,
 \end{itemize}
 such that
 \begin{enumerate}\itemsep=0pt
 \item[$1)$] $f_\alpha^\alpha=\id_{A_\alpha}$ for all $\alpha\in I$,
 \item[$2)$] $f^{\beta}_\gamma\circ f^{\alpha}_\beta=f^{\alpha}_\gamma$ for all $\alpha\leq\beta\leq\gamma$.
 \end{enumerate}
 An object $\mathscr A$ together with morphisms $f^\alpha\colon A_{\alpha}\to\mathscr A$ for $\alpha\in I$ is called \emph{inductive limit} of the inductive system $((A_\alpha)_{\alpha\in I},(f^{\alpha}_{\beta})_{\alpha\leq\beta})$ if
 \begin{enumerate}\itemsep=0pt
 \item[$1)$] $f^\alpha=f^{\beta}\circ f_{\beta}^\alpha$ for all $\alpha\leq \beta$,
 \item[$2)$] whenever $g^\alpha=g^{\beta}\circ f_{\beta}^\alpha$ holds for a family of morphisms $g^{\alpha}\colon A_{\alpha}\to B$ to some common object $B$, there exists a unique morphism $g\colon \mathscr A\to B$ such that $g\circ f^{\alpha}=g^{\alpha}$ for all $\alpha\in I$. This is referred to as the \emph{universal property} of the inductive limit.
 \end{enumerate}
\end{Definition}

If an inductive limit exists, it is essentially unique. More precisely, if $(\mathscr A,(f^{\alpha})_{\alpha\in I})$ and $(\mathscr B,(g^{\alpha})_{\alpha\in I})$ are two inductive limits of the same inductive system $(A_\alpha)_{\alpha\in I}$, then the uniquely determined morphisms $f\colon \mathscr A\to\mathscr B$ with $f\circ f^{\alpha}=g^{\alpha}$ and $g\colon \mathscr B\to\mathscr A$ with $g\circ g^{\alpha}=f^{\alpha}$ are mutually inverse isomorphisms.

In general, inductive limits may or may not exist. We call a category in which all inductive systems have inductive limits \emph{inductively complete}. See \cite[Chapter 12]{AHS04} for general arguments showing that many categories we consider even fulfill the stronger property of \emph{cocompleteness}.

A subset $J$ of a directed set $I$ is called \emph{cofinal} if for every $\alpha\in I$ there exists a $\beta\in J$ with $\beta\geq\alpha$.

\begin{Example}\label{sec:inductive-limits-1:ex:cofinal}
 For any fixed $\alpha_0\in I$ the set
 $\{\beta\mid\beta\geq\alpha_0\}$ is cofinal. Indeed, since $I$ is
 directed, there is a $\beta\geq\alpha,\alpha_0$ for all $\alpha\in
 I$.
\end{Example}

Clearly, if $\big((A_\alpha)_{\alpha\in I},(f^\alpha_\beta)_{\alpha\leq\beta,\alpha,\beta\in I}\big)$ is an inductive system and $J\subset I$ cofinal, then also $\big((A_\alpha)_{\alpha\in J},(f^\alpha_\beta)_{\alpha\leq\beta,\alpha,\beta\in J}\big)$ is an inductive system. It is known that the inductive limits are canonically isomorphic if they exist. We will need the following generalization of this (from which the mentioned fact follows as a special case as documented in \Cref{sec:inductive-limits-2:cor:cofinal} below).
Let $(A_\alpha)_{\alpha\in I}$ be an inductive system with inductive limit $(\mathscr A,(f^{\alpha})_{\alpha\in I})$, $K$ a directed set and $J_k\subset I$ for each $k\in K$ such that
 \begin{itemize}\itemsep=0pt
 \item $J_k$ is directed for all $k\in K$,
 \item $J_k\subset J_{k'}$ for all $k\leq k'$,
 \item $J=\bigcup_{k\in K} J_k$ cofinal in $I$.
 \end{itemize}
 Suppose the inductive systems $(A_\alpha)_{\alpha\in J_k}$ have inductive limits $\big(\mathscr A_k,(f^{\alpha}_{(k)})_{\alpha\in J_k}\big)$. It holds that $f^\alpha=f^\beta\circ f ^{\alpha}_{\beta}$ for all $\alpha\leq\beta\in J_k$, since $J_k\subset I$. Similarly, for $k\leq k'$ it holds that $f^\alpha_{(k')}=f^\beta_{(k')}\circ f ^{\alpha}_{\beta}$ for all $\alpha\leq\beta\in J_k$, since $J_k\subset J_{k'}$.
By the universal property of the inductive limit $\mathscr A_k$ there are unique morphisms $f^k_{k'}\colon \mathscr A_k\to\mathscr A_{k'}$ for $k\leq k'$ and $f^k\colon \mathscr A_k\to\mathscr A$ such that the diagrams
\[
\begin{tikzcd}[column sep=large]
 A_\alpha\rar{{f^\alpha}}\dar{{f^\alpha_{(k)}}}&\mathscr A, & A_\alpha\rar{{f^\alpha_{(k')}}}\dar{{f^\alpha_{(k)}}}&\mathscr A_{k'} \\
 \mathscr A_k\urar[swap]{{f^k}} & & \mathscr A_k\urar[swap]{{f^k_{k'}}}&
\end{tikzcd}
\]
commute for all $\alpha\in J_k$.

\begin{Proposition}\label{prop:ind-lim-agree}
 In the described situation $\big((\mathscr A_k)_{k\in K},(f^k_{k'})_{k\leq k'}\big)$ is an inductive system with inductive limit $\big(\mathscr A,(f^k)_{k\in K}\big)$.
\end{Proposition}

\begin{proof}
 The diagrams
\[
\begin{tikzcd}[column sep=large]
 A_\alpha\rar{{f^\alpha}}\dar{{f^\alpha_{(k)}}}\drar{f^{\alpha}_{(k')}}&\mathscr A
%%%%%%%%%%%%%%%%%
&
%%%%%%%%%%%%%%%%%
 A_\alpha\rar{{f^\alpha_{(k'')}}}\dar{{f^\alpha_{(k)}}}\drar{f^{\alpha}_{(k')}}&\mathscr A_{k''} \\
 \mathscr A_k \rar{f^k_{k'}} &\mathscr A_{k'},\uar[swap]{f^{k'}}
%%%%%%%%%%%%%%%%%%
&
%%%%%%%%%%%%%%%%%%%
\mathscr A_k\rar{f^k_{k'}}&\mathscr A_{k'}\uar[swap]{f^{k'}_{k''}}
\end{tikzcd}
\]
commute, which implies $f^{k}=f^{k'}\circ f^{k}_{k'}$ and $f^{k}_{k''}=f^{k'}_{k''}\circ f^{k}_{k'}$ for all $k\leq k'\leq k''$. Now suppose there are $g^k\colon A_k\to B$ with $g^k=g^{k'}\circ f^k_{k'}$ for all $k\leq k'$. We put $g^{\alpha}:= g^{k}\circ f^{\beta}_{(k)}\circ f^{\alpha}_\beta$ for $\beta\in J_k$, $\alpha\leq\beta$. Since $J=\bigcup_{k\in K} J_k$ is cofinal in $I$, we can find such $k$ and $\beta$ for every $\alpha$. One can check that the $g^{\alpha}$ do not depend on the choice and fulfill $g^{\alpha}=g^{\alpha'}\circ f^{\alpha}_{\alpha'}$ for all $\alpha\leq\alpha'$. This yields a~morphism $g\colon \mathscr A\to\mathscr B$ which makes
\[
\begin{tikzcd}[column sep=large]
 \mathscr A_k\rar{g^{k}}\drar{f^{k}}&\mathscr B \\
 A_\beta \rar{f^\beta}\uar{{f^\beta_{(k)}}} &\mathscr A\uar[swap]{g}\\
 A_\alpha\uar{f^{\alpha}_{\beta}}\urar[swap]{f^{\alpha}}&
\end{tikzcd}
\]
commute. On the other hand, any morphism which makes the upper right triangle commute, automatically makes the whole diagram commute and will therefore equal $g$.
\end{proof}

\begin{Corollary}\label{sec:inductive-limits-2:cor:cofinal}
 Let $(A_\alpha)_{\alpha\in I}$ be an inductive system with inductive limit $\big(\mathscr A,(f^{\alpha})_{\alpha\in I}\big)$, $J\subset I$ cofinal. Then $(A_\alpha)_{\alpha\in J}$ is an inductive system with inductive limit $\big(\mathscr A,(f^{\alpha})_{\alpha\in J}\big)$.
\end{Corollary}

\begin{proof}
 This is a special case of \Cref{prop:ind-lim-agree} with $\betrag K=1$, since the inductive system over the one point set $K$ does not add anything.
\end{proof}

\subsection{Tensor categories}\label{sec:categ-indep}

A \emph{tensor category} is a category $\mathcal C$ together with a bifunctor $\boxtimes\colon \mathcal C\times\mathcal C\to\mathcal C$ which
\begin{itemize}\itemsep=0pt
\item is associative under a natural isomorphism with components
\[
\alpha_{A,B,C}\colon \ A\boxtimes (B\boxtimes C)\xrightarrow[]{\cong} (A\boxtimes B)\boxtimes C
\] called \emph{associativity constraint},
\item has a unit object $E\in\Obj(\mathcal C)$ acting as left and right identity under natural isomorphisms with components
\[
l_A\colon \ E\boxtimes A\xrightarrow{\cong} A,\qquad r_A\colon \ A\boxtimes E\xrightarrow{\cong} A
\]
called \emph{left unit constraint} and \emph{right unit constraint} respectively
\end{itemize}
such that the pentagon and triangle identities hold~\cite{McL98}. If the natural transformations $\alpha$, $l$ and~$r$ are all identities, we say the tensor category is \emph{strict}.

It can be shown that the pentagon and triangle identities imply commutativity of all diagrams which only contain~$\alpha$,~$l$ and~$r$ \cite[Section~VII.2]{McL98}. This is called Mac Lane's \emph{coherence} theorem. Even for non-strict tensor categories, we will frequently suppress the associativity and unit constraints in the notation and write $(\mathcal C,\boxtimes,E)$, or even $(\mathcal C,\boxtimes)$ or $\mathcal C$. In the examples we treat, $\alpha$, $l$ and $r$ are always canonical.

Given tensor categories $(\mathcal C,\boxtimes)$ and $(\mathcal C',\boxtimes')$ with unit objects, associativity and unit constraints $E$, $\alpha$, $l$, $r$ and $E'$, $\alpha'$, $l'$, $r'$ respectively,
a \emph{cotensor functor} is a triple $(\mathcal F,\delta,\Delta)$ consisting of
\begin{itemize}\itemsep=0pt
\item a functor $\mathcal F\colon \mathcal C\to\mathcal C'$,
\item a morphism $\delta\colon \mathcal F(E)\to E'$,
\item a natural transformation $\Delta\colon \mathcal F(\cdot\,\boxtimes\,\cdot)\Rightarrow\mathcal F(\cdot)\boxtimes'\mathcal F(\cdot)$
\end{itemize}
such that the diagrams
\begin{equation}\label{GLa15:eq:cotensor1}
\begin{tikzcd}[column sep=7.2em]
\mathcal F\bigl(A\boxtimes(B\boxtimes C)\bigr) \rar{\mathcal F(\alpha_{A,B,C})}\dar[swap]{\Delta_{A,B\boxtimes C}}
&\mathcal F\bigl((A\boxtimes B)\boxtimes C\bigr)\dar{\Delta_{A\boxtimes B,C}}\\
\mathcal F(A)\boxtimes'\mathcal F(B\boxtimes C)\dar[swap]{\id_{\mathcal F(A)}\boxtimes'\Delta_{B,C}}
& \mathcal F(A\boxtimes B)\boxtimes'\mathcal F(C)\dar{\Delta_{A,B}\boxtimes'\id_{\mathcal F(C)}}\\
\mathcal F(A)\boxtimes'\bigl(\mathcal F(B)\boxtimes'\mathcal F(C)\bigr)
\rar{\alpha'_{\mathcal F(A),\mathcal F(B),\mathcal F(C)}}
&\bigl(\mathcal F(A)\boxtimes'\mathcal F(B)\bigr)\boxtimes'\mathcal F(C),
\end{tikzcd}
\end{equation}
\begin{equation}\label{GLa15:eq:cotensor2}
\begin{tikzcd}
 \mathcal F(B\boxtimes E)\rar{\Delta_{B,E}}\dar[swap]{\mathcal
 F(r_B)} & \mathcal F(B)\boxtimes'\mathcal
 F(E)\dar{\id_{\mathcal F(B)}\boxtimes' \delta}
 &%%%%%%%%%%%%%%%%%%
 \mathcal F(E\boxtimes B)\rar{\Delta_{E,B}}\dar[swap]{\mathcal
 F(l_B)} & \mathcal F(E)\boxtimes'\mathcal
 F(B)\dar{\delta\,\boxtimes'\id_{\mathcal F(B)}}
 \\
 \mathcal F(B) & \mathcal F(B)\boxtimes' E',\lar[swap]{r'_{\mathcal
 F(B)}} &%%%%%%%%%%%%%%%%%%
 \mathcal F(B) & E'\boxtimes'\mathcal F(B)\lar[swap]{l'_{\mathcal
 F(B)}}
\end{tikzcd}
\end{equation}
commute for all $A,B,C\in\Obj(\mathcal C)$. A cotensor functor is called \emph{strong} if $\Delta$ is a natural isomorphism and $\delta$ is an isomorphism.

\begin{Theorem}\label{sec:categ-indep-1:cotensor-functor-composition}
 Let $\mathcal F\colon \mathcal C\to\mathcal C'$ and $\mathcal F'\colon \mathcal C'\to\mathcal C''$ be cotensor functors with coproduct morphisms $\Delta_{A,B}$, $\Delta'_{A',B'}$ and counit morphisms $\delta$, $\delta'$. Then $\mathcal F'\circ\mathcal F$ is a cotensor functor with coproduct morphisms $\Delta'_{\mathcal F(A),\mathcal F(B)}\circ\mathcal F'(\Delta_{A,B})$ and counit morphism $\delta'\circ\mathcal F'(\delta)$.
\end{Theorem}

This is well known and can be shown by writing down the involved diagrams and check that they commute; see~\cite{Lac15} for an explicit proof.

Similarly, a \emph{tensor functor} is a functor $\mathcal F\colon \mathcal C\to\mathcal C'$ together with a natural transformation $\mu\colon \mathcal F(\cdot)\boxtimes'\mathcal F(\cdot)\Rightarrow\mathcal F(\cdot\boxtimes \cdot)$ and a morphism $\Eins\colon E'\to\mathcal F(E)$ such that the diagrams one obtains from~\eqref{GLa15:eq:cotensor1} and~\eqref{GLa15:eq:cotensor2} by reversing the arrows and replacing $\Delta$ and $\delta$ with $\mu$ and $\Eins$ commute.

\section{Independence in tensor categories}\label{sec:indep-tens-categ}

In this section we quickly recall Franz' original definition of independence in tensor categories with inclusions (Section~\ref{sec:tensor-categories-1}), prove equivalence of different compatibility conditions between the inclusions and initiality of the unit object (Section~\ref{sec:comp-incl}), and explore what compatibility means for the obvious extension of Franz' independence to more than two morphisms (Section~\ref{sec:good-inclusions}).

\subsection{Categorial independence with respect to inclusions}\label{sec:tensor-categories-1}

In order to unify the different notions of independence in quantum probability, Franz came up with a definition of independence in a tensor-categorial framework \cite[Section~3]{Fra06}. Let $\mathcal P_i\colon \mathcal C\times\mathcal C\allowbreak \to\mathcal C$ for $i\in\{1,2\}$ denote the projection functor onto the first or second component respectively.

\begin{Definition}\label{def:tensorcat_w_inclusions}
Let $(\mathcal C,\boxtimes)$ be a tensor category. A natural transformation
 $\iota^{1}\colon \mathcal P_{1}\Rightarrow{\boxtimes}$ is called \emph{left inclusion} and a natural transformations
 $\iota^{2}\colon \mathcal P_{2}\Rightarrow{\boxtimes}$ is called \emph{right inclusion}. A~tensor category together with a right and a left inclusion is referred to as \emph{tensor category with inclusions}.
\end{Definition}

In more detail, inclusions for a tensor category are two collections of morphisms $\iota^i_{B_1,B_2}\colon B_{i}\to B_1\boxtimes B_2$ for $B_1,B_2\in\mathcal C$, $i\in\{1,2\}$ such that
\[
\begin{tikzcd}[column sep=huge]
A_1\dar{f_1}\rar{\iota^1} & A_1\boxtimes A_2 \dar{f_1\boxtimes f_2}
& A_2\dar{f_2}\lar[swap]{\iota^2}\\
B_1 \rar{\iota^1}& B_1\boxtimes B_2 & B_2\lar[swap]{\iota^2}
\end{tikzcd}
\]
commutes for all $f_{i}\colon A_{i}\to B_{i}$, $i\in\{1,2\}$.

\begin{Definition}\label{sec:categ-indep-1:def:cat-ind}
 Let $\big(\mathcal C,\boxtimes,\iota^1,\iota^2\big)$ be a tensor category with inclusions. Two morphisms $j_1,j_2$: $B_i\to A$ are \emph{independent} if there exists a morphism $h\colon B_1\boxtimes B_2\to A$ such that the diagram
\[
\begin{tikzcd}[column sep=huge]
 \ & A &\ \\
 B_1\urar{j_1}\rar{\iota^{1}} & B_1\boxtimes B_2\uar[swap]{h}
 & B_2\ular[swap]{j_2}\lar[swap]{\iota^{2}}
\end{tikzcd}
\]
commutes. Such a morphism $h$ is called \emph{independence morphism} for $j_1$ and $j_2$.
\end{Definition}

% \begin{Example}
% We begin with a trivial example. The direct sum of vector spaces $\mathcal V_1$ and $\mathcal V_2$ has the following property: Given any two linear maps $f_i\colon \mathcal V_i\to\mathcal W$ to some third vector space $\mathcal W$, there exists a unique linear map $h\colon \mathcal V_1\oplus\mathcal V_2\to\mathcal W$ with $h(v_i)=f_i(v_i)$ for all $v_i\in\mathcal V_i$, namely $h=f_1+f_2$ (here we identify $V_i$ with the corresponding subspace of $\mathcal V_1\oplus\mathcal V_2$). In particular, in the tensor category $(\Vect,\oplus)$ with the canonical inclusions $\mathcal V_i\hookrightarrow\mathcal V_1\oplus\mathcal V_2$, all pairs of linear maps into a common vector space are independent.

% In the tensor category $(\Alg,\sqcup)$ with $\sqcup$ the free product of algebras the situation is similar. With respect to the canonical inclusions $\mathcal A_i\hookrightarrow\mathcal A_1\sqcup\mathcal A_2$, any pair of algebra homomorphisms $j_i\colon \mathcal A_i\to\mathcal B$ to a common algebra $\mathcal B$ is independent.
% \end{Example}

Let us first consider the case where the tensor product $B_1\boxtimes B_2$ coincides with the \emph{coproduct} $B_1\sqcup B_2$ in the category; see \cite{Fra06} or \cite{McL92}. By definition of a coproduct, any pair of morphisms $j_i\colon B_i\to A$ to a common target $A$ will be independent with the unique independence morphism $h=j_1\sqcup j_2$. Coproducts exist in many categories, for example the direct sum in the category of vector spaces with linear maps, or the free product in the category of algebras with algebra homomorphisms. In order to have a nontrivial notion of independence, one should either use a~different tensor product, or restrict the class of morphisms. In \Cref{sec:examples-1} we will see that many notions of independence used in mathematics are indeed special cases of categorial independence, in particular this holds for linear independence, orthogonality, stochastic independence, Bose or tensor independence, Fermi independence, and all notions of noncommutative stochastic independence which are induced by universal products, like freeness, Boolean independence and monotone independence.

In most examples the independence morphism $h$ will be uniquely determined if it exists. The next example shows that this is not the case in general, which is why we will not assume uniqueness in the development of the general theory.

\begin{Example} Consider the category $\Vect$ with tensor product
 \begin{align*}
 V_1 \odot V_2:= V_1\oplus V_2\oplus V_1\otimes V_2
 \end{align*}
 and the canonical inclusions $V_i\hookrightarrow V_1\odot V_2$ which identify $V_i$ with the summand $V_i$ in $V_1\oplus V_2$. Any two linear maps $f_i\colon V_i\to W$ are independent, but the independence morphism is not uniquely determined. Indeed, for an arbitrary linear map $f\colon V_1\otimes V_2\to W$, the linear map $h=f_1+f_2+f$ is an independence morphism for $f_1$ and $f_2$.
\end{Example}

\subsection{Compatible inclusions}\label{sec:comp-incl}

We already defined what it means for two morphisms $j_1$, $j_2$ in a tensor category with inclusions to be independent. But if we want to consider notions of independence for more than two morphisms, we need to require certain compatibility conditions between the inclusions and the structure of the tensor category. In particular, for dealing with categorial L\'evy processes (as we will do in Section~\ref{sec:second-induct-limit}) it seems to be necessary. This has first been observed in~\cite{G15} and the conditions do not appear in~\cite{Fra06}.

\begin{Definition}
 Inclusions $\iota^1$, $\iota^2$ are called \emph{compatible with the unit constraints} if the diagram\looseness=-1
 \begin{equation}
 \begin{tikzcd}[column sep=huge]
 E\boxtimes A\arrow{dr}[swap]{l_A} &A\lar[swap]{\iota^2}\rar{\iota^1}\dar{\id_A}&A\boxtimes E\arrow{dl}{r_A}\\
 &A&
 \end{tikzcd}
 \label{eq:compatible-w-unit-constraints}
\end{equation}
 commutes for all objects $A\in\mathcal{C}$.
\end{Definition}

\begin{Theorem}\label{theo:compincl-uoinitial}
 Let $\mathcal C$ be a tensor category.
 \begin{enumerate}\itemsep=0pt
 \item[$(a)$] %\label{item:comp-incl-th-a}
 If $\iota^1$, $\iota^2$ are inclusions which are compatible with the unit constraints, then the unit object $E$ is initial, i.e., there is a unique morphism $\Eins_A\colon E\to A$ for every object $A\in\mathcal{C}$. Furthermore,
 \begin{align}\label{eq:incl-from-init}
 \iota^1_{A,B}=(\id\boxtimes\Eins_B)\circ r_A^{-1},
 \qquad
 \iota^2_{A,B}=(\Eins_A\boxtimes\id)\circ l_B^{-1}
 \end{align}
 holds for all objects $A,B\in\mathcal{C}$.
 \item[$(b)$] %\label{item:comp-incl-th-b}
 Suppose that the unit object $E$ is an initial object. Then~\eqref{eq:incl-from-init} read as a definition yields inclusions $\iota^1$, $\iota^2$ which are compatible with the unit constraints.
\end{enumerate}
\end{Theorem}

\begin{Remark} After writing this article, we learned from Tobias Fritz and an anonymous referee that the discussed relationship between inclusions and an initial unit object (or, dually, projections and a terminal unit object) have been known to many category theorists, but we do not know of any reference where it is spelled out (cf.\ Tobias Fritz' post \cite{website-golem-semicartesian} on \emph{golem} and the \emph{nLab} page \cite{website-nLab-semicartesian} on semicartesian categories). We therefore think of it as a valuable contribution to fill this gap.
\end{Remark}

\begin{proof}(a) The inclusions can easily be used to define the morphism $\Eins_A:=l_A\circ\iota^1_{E,A}$ from $E$ to an arbitrary object $A$. To prove that $E$ is initial, it remains to show that any morphism $f\colon E\to A$ coincides with $\Eins_A$. Naturality of $l$ yields $f\circ l_E=l_A\circ (\id\boxtimes f)$. The unit constraints $r_E,l_E\colon E\otimes E\to E$ coincide in any tensor category by coherence. From the compatibility with the unit constraints, \eqref{eq:compatible-w-unit-constraints} for $A=E$, we conclude that $\iota^1_{E,E}=\iota^2_E\colon E\to E\boxtimes E$ and $r_E=l_E\colon E\boxtimes E \to E$ are mutually inverse isomorphisms. Thus, solving for $f$ yields $f=l_A\circ (\id\boxtimes f)\circ \iota^1_{E,E}$. Naturality of $\iota^1$ implies that $(\id\boxtimes f)\circ\iota^1_{E,E}=\iota^1_{E,A}$. Therefore $f=l_A\circ\iota^1_{E,A}=\Eins_A$, which shows that $E$ is initial.
In the diagram
\[
\begin{tikzcd}[column sep=huge]
 A\arrow{rr}{\iota^1}\drar{\iota^1}\arrow{dddr}[swap]{r^{-1}=\iota^1}&&A\boxtimes B\\
 &A\boxtimes(B\boxtimes E)\urar{\id\boxtimes r}&\\
 &&\\
 &A\boxtimes E\arrow{uu}{\id\boxtimes\iota^2}\arrow{uuur}[swap]{\id\boxtimes\Eins_B}&
\end{tikzcd}
\]
the upper triangle and the lower left triangle commute due to the naturality of $\iota^1$. Because $E$ is initial, it holds that $\Eins_B=r_A\circ\iota^2_{A,E}$, so the lower right triangle also commutes. So the whole diagram commutes and the outside triangle represents the first equation of \eqref{eq:incl-from-init}. The second one follows analogously.

(b) It follows from uniqueness of morphisms from $E$ to $B$ that $\Eins_B=\Eins_A\circ f$ for all $f\colon A\to B$, or put as a diagram,
\[
\begin{tikzcd}
 E\dar[swap]{\id_E}\rar{\Eins_A}&A\dar{f}\\E\rar[swap]{\Eins_B}&B
\end{tikzcd}
\]
commutes. Thus, we can interpret the collection of all $\Eins_A,A\in\mathcal{C}$ as a natural transformation $\Eins\colon E\Rightarrow \id_\mathcal{C}$; here $E$ stands for the constant functor $E\colon\mathcal{C}\to\mathcal{C}$ with $E(A)=E$ for all objects $A$ and $E(f)=\id_E$ for all morphisms $f$. With this, naturality of $\iota^1,\iota^2$ is easy to check. We also have that $\Eins_E=\id_E$. Hence, $r_A\circ\iota_{A,E}=r_A\circ(\id_A\boxtimes\Eins_E)\circ r_A^{-1}=r_Ar_A^{-1}=\id_A$ for all objects $A$. Analogously, $l_A\circ\iota_{E,A}^2=\id_A$.
\end{proof}

\begin{Observation}\label{obs:coherence}
 Mac Lanes's coherence theorem extends to diagrams built up from the natural transformations $\alpha$, $l$, $r$ and $\Eins$. For a formal statement, there are some subtleties to consider, mainly caused by the problem that in a concrete tensor category it can happen that tensor products of different factors can yield the same object. Also, one has to be very careful what is meant by ``built up from natural transformations''; e.g., $\iota^1_{A,A},\iota^2_{A,A}\colon A\to A\boxtimes A$ will typically not be equal.

 Concretely, we will need commutativity of the diagrams
\begin{equation}
 \begin{tikzcd}[column sep=huge]
 A\boxtimes C\arrow{d}[swap]{\id_A\boxtimes\iota^2}\rar{\iota^1\boxtimes\id_C}&(A\boxtimes B)\boxtimes C,\\ A\boxtimes(B\boxtimes C)\arrow{ur}[swap]{\alpha_{A,B,C}}&
 \end{tikzcd}
 \label{eq:compatible-w-asso-constraints}
\end{equation}
 \[
 \begin{tikzcd}[column sep=small]
 A\boxtimes B\dar{\iota^1}&B\rar{\iota^1}\lar[swap]{\iota^2}&B\boxtimes C\dar{\iota^2}
 &%%%%%%%%%%%%%%%%%%
C\rar{\iota^2}\dar{\iota^2}&B\boxtimes C\dar{\iota^2}\\
 (A\boxtimes B)\boxtimes C\arrow{rr}{\alpha}&&A\boxtimes (B\boxtimes C),
 &%%%%%%%%%%%%%%%%%%%%%%%%%%%%%%%
 (A\boxtimes B)\boxtimes C\rar{\alpha}&A\boxtimes( B\boxtimes C),
 \end{tikzcd}
\]
more generally, for $\{i_1<\dots< i_k\}\subset\{1,\dots, n\}$ $($suppressing associativity constraints$)$ there is a unique natural transformation $\iota^{i_1,\dots, i_k;n}$, called inclusion $($and geralizing $\iota^1=\iota^{1;2},\iota^2=\iota^{2;2})$,
\[\iota^{i_1,\dots, i_k;n}_{A_1,\dots, A_n}\colon \ A_{i_1}\boxtimes\cdots\boxtimes A_{i_k}\to A_1\boxtimes\cdots\boxtimes A_n, \qquad A_1,\dots, A_n\in\mathcal C,\]
which can be written as $($vertical\footnote{The vertical composition is just the componentwise composition, i.e., if $\alpha\colon F_1\Rightarrow F_2$ and $\beta\colon F_2\Rightarrow F_3$, then $\beta\circ\alpha\colon F_1\Rightarrow F_3$ has components $\beta_{A}\circ \alpha_A\colon F_1(A)\to F_3(A)$}$)$ composition of tensor products with $\mathrm{id}$ of
\begin{itemize}\itemsep=0pt
\item $l$, $r$, $\Eins$ $($and $\alpha$ if we remembered the parentheses$)$,
\item the ``degenerate'' natural transformations obtained by replacing an argument by the unit object $E$, e.g., $\alpha_{\cdot,E,\cdot}=(\alpha_{A,E,C})_{A,C\in\mathcal C}$ or $r_E$,
\item the inverses of all those natural transformations except $\Eins$ $($which is the only one that is not a natural isomorphism$)$.
\end{itemize}
The basic idea to prove our claim is the following. Suppose that $f_1$ and $f_2$ are parallel natural transformations of the described form. We can commute instances of $\Eins$ with instances of $l$, $r$ $($and~$\alpha)$ using naturality, for example $r\circ(\Eins\boxtimes\id_E)=\Eins\circ r_E$ or $\alpha\circ ((\id\boxtimes\Eins)\boxtimes \id)=(\id\boxtimes(\Eins\boxtimes \id))\circ\alpha_{\cdot,E,\cdot}$. This allows us to subsequently move all instances of $\Eins$ to the left of all instances of~$\alpha$, $l$ and $r$. Thus, we can factorize $f_i$ as $h_i\circ g_i$ with $g_i$ built up from~$\alpha$, $l$ and $r$ whereas $h_i$ is built up from~$\Eins$ alone. Necessarily, with the notation
\[\tilde A_j:=
 \begin{cases}
 A_j, & j\in \{i_1,\dots, i_k\},\\
 E, & j\notin \{i_1,\dots, i_k\},\end{cases}
 \qquad
 \tilde\id_j=\begin{cases}
 \mathrm{id},& j\in \{i_1,\dots, i_k\},\\
 \Eins, & j\notin \{i_1,\dots, i_k\}
 \end{cases}\]
 we have
\begin{align*}
 h_1=h_2= \tilde\id_1 \boxtimes \dots \boxtimes \tilde\id_n,\qquad (h_i)_{A_1,\dots, A_n}\colon \ \tilde A_1\boxtimes\cdots\boxtimes \tilde A_n\to A_1\boxtimes\dots\boxtimes A_n.
\end{align*}
It follows that $g_1$ and $g_2$ are parallel,
\[(g_i)_{A_1,\dots, A_n}\colon \ A_{i_1}\boxtimes\cdots\boxtimes A_{i_k}\to \tilde A_1\boxtimes\cdots\boxtimes \tilde A_n,
\]
so they are equal by Mac Lane's coherence theorem. Therefore $f_1=h_i\circ g_i=f_2$.
\end{Observation}

\subsection{General categorial independence}\label{sec:good-inclusions}

Let $\mathcal{C}$ be a tensor category such that the unit object $E$ is an initial object. By \Cref{obs:coherence} there are unique natural inclusions $\iota^{i_1,\dots, i_k;n}_{A_1,\dots, A_n}\colon A_{i_1}\boxtimes\cdots\boxtimes A_{i_k}\to A_1\boxtimes\cdots\boxtimes A_n$ built up from $l$, $r$ and~$\Eins$.

\begin{Definition}
 Let $B_1,\dots, B_n,A$ be objects of $\mathcal{C}$ and $f_i\colon B_i\to A$ morphisms. Then $f_1,\dots,f_n$ are called \emph{independent} if there exists a morphism $h\colon B_1\boxtimes\cdots\boxtimes B_n\to A$ such that the diagrams
\[
\begin{tikzcd}
 \ & A \\
 B_i\rar[swap]{\iota^{i;n}}\arrow{ur}{f_i}&B_1\boxtimes\cdots\boxtimes B_n\uar{h}
\end{tikzcd}
\]
commute for all $i\in\{1,\dots, n\}$; $h$ is then called an \emph{independence morphism} for $f_1,\dots, f_n$.
\end{Definition}

We conclude with the analogues of two basic results about stochastic independence: Subfamilies of independent families of random variables are independent and functions of independent random variables are independent.

\begin{Theorem} Let $\mathcal{C}$ be a tensor category with initial unit object. Furthermore let $f_1,\dots, f_n$, $f_i\colon B_i\to A$, be independent with independence morphism $h\colon B_1\boxtimes\cdots\boxtimes B_n\to A$. Then the following holds.
 \begin{enumerate}\itemsep=0pt
 \item[$(a)$] %\label{item:subfam-ind}
 For all $1\leq i_1<\dots<i_k\leq n$ the morphisms $f_{i_1},\dots,f_{i_n}$ are independent with independence
 morphism $h\circ\iota^{i_1,\dots,i_k;n}$.
 \item[$(b)$] %\label{item:func-ind}
 Let $j_1,\dots,j_1$, $j_i\colon C_i\to B_i$, be morphisms and put $g_i:=f_i\circ j_i\colon C_i\to A$. Then $g_1,\dots, g_n$ are independent with independence morphism $h\circ (j_1\boxtimes\cdots\boxtimes j_n)$.
 \end{enumerate}
\end{Theorem}

\begin{proof}
(a) In the diagram
 \[
 \begin{tikzcd}[column sep=huge]
 \ &A\\
 B_{i_{\!j}}\arrow{ur}{f_{i_{\!j}}}\rar{\iota^{i_{\!j};n}}\arrow{dr}{\iota^{j;k}}&B_1\boxtimes\cdots\boxtimes B_n\uar[swap]{h}\\
 \ &B_{i_1}\boxtimes\cdots \boxtimes B_{i_k}\uar[swap]{\iota^{i_1,\dots,i_k;n}}
\end{tikzcd}
\]
the upper half commutes by independence of $f_1,\dots, f_n$ and the lower half by coherence (see Observation \ref{obs:coherence}). Commutativity for every $j=1,\dots,k$ proves the assertion.

(b) In the diagram
 \[
 \begin{tikzcd}[column sep=huge]
 C_i\rar{j_i}\dar{\iota^{i;n}} &B_i\dar{\iota^{i;n}}\arrow{dr}{f_i}&\ \\
 C_1\boxtimes\cdots\boxtimes C_n\rar{j_1\boxtimes\cdots\boxtimes j_n}&B_1\boxtimes\cdots\boxtimes B_n\rar{h}&A
\end{tikzcd}
\]
the left hand side square commutes because $\iota^{i;n}$ is a natural transformation and the right hand side triangle commutes by independence of $f_1,\cdots,f_n$. So we get $g_i=f_i\circ j_i= h\circ (j_1\boxtimes\cdots\boxtimes j_n)\circ \iota^{i;n}_{C_1,\dots,C_n}$ for all $i$, which was the assertion.
\end{proof}

\begin{Remark}
 We formulated (a) and (b) the way we will use them later. However, note that under the canonical identification of tensor products with the initial unit object, (a) can also be seen as a special case of (b) with $j_i=
 \begin{cases}
 \id_{B_{i_k}}, & \text{$i=i_k$ for some $k\in\{1,\dots,n\}$},\\
 \Eins_A, & \text{otherwise.}
 \end{cases}$

 Although we will not make use of it in the sequel, it is also noteworthy that (b) can be generalized as follows, thus relating it to the main condition in Simpson's definition of independence structure \cite[Definition 2.1]{Sim18}:
 \begin{enumerate}\itemsep=0pt
 \item[(c)] Suppose that
 \begin{itemize}\itemsep=0pt
 \item $f_k\colon B_k\to A$, $k\in\{1,\dots, n\}$ are independent with independence morphism $h$,
 \item $f_{k,\ell}\colon C_{k,\ell}\to B_k$, $\ell\in\{1,\dots, m_k\}$ are independent with independence morphism $j_k\colon C_k:=C_{k,1}\boxtimes \dots\boxtimes C_{k,m_k}\to B_k$ for each $k\in\{1,\dots, n\}$.
 \end{itemize}
 Then $g_{k,\ell}:= f_k\circ f_{k,\ell}\colon C_{k,\ell}\to A$ are independent with independence morphism
 \[h\circ (j_1\boxtimes \dots\boxtimes j_n).\]
 \end{enumerate}
 Indeed, in the diagram
 \[
 \begin{tikzcd}[column sep=huge]
 C_{k,\ell}\rar{f_{k,\ell}}\dar{\iota^{\ell;m_k}} &B_k\dar{\iota^{k;n}}\arrow{r}{f_k}&A \\
 C_{k}\arrow{ur}{j_k}\arrow{d}{\iota^{k;n}}&B_1\boxtimes\cdots\boxtimes B_n\arrow{ur}{h}&\ \\
 C_1\boxtimes\dots \boxtimes C_n\arrow{ur}{j_1\boxtimes\dots \boxtimes j_n}&\ &\
\end{tikzcd}
\]
the upper left triangle commutes by definition of independence morphisms and the rest of the diagram commutes as noted in the proof of (b).
\end{Remark}

We briefly discuss independence for infinite families. Let $(f_i\colon B_i\to A )_{i\in I}$ be a family of morphisms indexed by a totally ordered index set $(I,{\leq})$. The set $P_{fin}(I)=\bigl\{\{i_1<\dots<i_n\}\mid n\in\mathbb N,i_1,\dots,i_n\in I\bigr\}$ of all finite subsets of $I$ is a directed set with respect to inclusion. For $\mathbf{i}=\{i_1<\cdots<i_n\}\in P_{fin}(I)$ $B_{\mathbf{i}}:=B_{i_1\boxtimes\cdots\boxtimes i_n}$ and $f_{\mathbf{i}}:=f_{i_1}\boxtimes\cdots\boxtimes f_{i_n}\colon B_{\mathbf{i}}\to A$. Given two finite subsets $\mathbf{i},\mathbf{j}\in P_{fin}(I)$ with $\mathbf{i}\subset\mathbf{j}$, we put $\iota_{\mathbf{i}}^{\mathbf{j}}\colon B_{\mathbf{i}}\to B_{\mathbf{j}}$ as the unique morphism described in \Cref{obs:coherence}. By the same observation we know that $\iota_{\mathbf{j}}^{\mathbf{k}}\iota_{\mathbf{i}}^{\mathbf{j}}=\iota_{\mathbf{i}}^{\mathbf{k}}$ for all $\mathbf{i}\subset \mathbf{j}\subset \mathbf{k}$. In other words, $\bigl((B_{\mathbf{i}})_{\mathbf{i}\in P_{fin}(I)},(f_{\mathbf{i}}^{j})_{\mathbf{i}\subset\mathbf{j}}\bigr)$ is an inductive system. We say that the $(f_i)_{i\in I}$ are independent if there is an inductive limit $\big(B_{I},(\iota_{\mathbf{i}\colon B\mathbf{i}\to B_I})_{\mathbf{i}\in P_{fin}(I)}\big)$ of $\bigl((B_{\mathbf{i}})_{\mathbf{i}\in P_{fin}(I)},(\iota_{\mathbf{i}}^{j})_{\mathbf{i}\subset\mathbf{j}}\bigr)$ and a morphism $h\colon B_I\to A$ such that $h\circ \iota^{\{i\}}=f_i$ for all $i\in I$.\footnote{Note that this inductive limit coincides with the \emph{infinite tensor product} in $\mathcal C^{\rm op}$ as defined in \cite{RiFr20} in case the tensor category $\mathcal C$ has the additional structure of a symmetric monoidal category.}

\section{L\'evy processes in tensor categories}\label{sec:comonoidal-systems}

In this section we define comonoidal systems and categorial L\'evy processes. We present two important inductive limit constructions which have been considered for many different examples and to which we give a unifying framework. In the end, these allow us to reconstruct under suitable conditions a categorial L\'evy process from its comonoidal system, which we like to think of as an analogue of the family of marignal distributions of a classical L\'evy process. In the context of Hilbert modules these constructions have first been described in full detail by Bhat and Skeide~\cite{BhSk00}, who also refer to them as first and second inductive limit.

We need a number of results on monoids, which we prove in Sections~\ref{sec:uf-monoids} and~\ref{sec:ore-monoids}. We do not claim originality of these results, however, we were also not able to find a suitable reference containing all the statements we need. Shalit and Skeide \cite[second half of Section~4, in particular Theorem~4.13]{ShSk20p} present a similar discussion, but it seems that the forward direction of our \Cref{thm:totally_ordered-iff-oreuf} was not known to them (see \Cref{rem:ShSk-question}).

We will conclude the section with a summary of the assumptions we make on the involved category and monoid in order to perform the different discussed constructions, see Section~\ref{sec:assumptions}.

\subsection{Comonoidal systems}\label{sec:comonoidal-systems-1}

A \emph{monoid} is a semigroup with a unit element. We identify a monoid~$\mathbb S$ with the strict tensor category whose objects are the elements of~$\mathbb S$ with only the identity morphisms and the tensor product given by the multiplication of~$\mathbb S$.

\begin{Definition}\label{sec:comonoidal-systems-1:Def:com-syst}
 Let $\mathbb S$ be a monoid and $(\mathcal C,\boxtimes)$ a tensor category. A \emph{monoidal system} over $\mathbb S$ in $\mathcal C$ is a tensor functor from $\mathbb S$ to $\mathcal C$. A \emph{comonoidal system} over $\mathbb S$ in $\mathcal C$ is a cotensor functor from $\mathbb S$ to $\mathcal C$. A comonoidal system is called \emph{full} if the cotensor functor is strong. A monoidal system (respectively comonoidal system) over the trivial monoid $\{e\}$ is simply called a \emph{monoid in $\mathcal C$} (respectively \emph{comonoid in $\mathcal C$}).
\end{Definition}

Since there are only identity morphisms in $\mathbb S$, any functor defined on $\mathbb S$ acts trivially on morphisms, so it is determined by the object assignment and can be identified with the family $(A_s)_{s\in\mathbb S}$ where $A_s$ denotes the value of the functor at $s\in\mathbb S$. Thus, a monoidal system over $\mathbb S$ in $\mathcal C$ is the same as a family of objects $(A_s)_{s\in\mathbb S}$ together with \emph{product morphisms} $\mu_{s,t}\colon A_s\boxtimes A_t\to A_{st}$ and a \emph{unit morphism} $u\colon E\to A_e$ such that the natural associativity and unit properties
\[
\begin{tikzcd}
 A_r\boxtimes A_s\boxtimes A_t\rar{\mu_{r,s}\boxtimes\id_t}\dar{\id_r\boxtimes\mu_{s,t}}&A_{rs}\boxtimes A_t\dar{\mu_{rs,t}}
%%%%%%%%%%%%%%%%
&
%%%%%%%%%%%%%%%%%%
E\boxtimes A_s\dar{u\boxtimes \id}&A_s\lar[swap]{l_{A_s}^{-1}}\rar{r_{A_s}^{-1}}\dar{\id}&A_s\boxtimes E\dar{\id\boxtimes u}
\\
 A_r\boxtimes A_{st}\rar{\mu_{r,st}}&A_{rst},
%%%%%%%%%%%%%%%%%%%
&
%%%%%%%%%%%%%%%%%%
A_e\boxtimes A_s\rar{\mu_{e,s}}&A_s&A_s\boxtimes A_e\lar[swap]{\mu_{s,e}}
\end{tikzcd}
\]
are fulfilled. % In particular, a monoid in $\Set$ is just a usual monoid.
Similarly, a comonoidal system over $\mathbb S$ in $\mathcal C$ is a family of objects $(A_s)_{s\in\mathbb S}$ together with \emph{coproduct morphisms} $\Delta_{s,t}\colon A_{st} \to A_s\boxtimes A_t$ and a \emph{counit morphism} $\delta\colon A_e\to E$ such that coassociativity and the counit properties
\[
\begin{tikzcd}
 A_{rst}\rar{ \Delta_{rs,t}}\dar{ \Delta_{r,st}}&A_{rs}\boxtimes A_t\dar{ \Delta_{r,s}\boxtimes\id_t}
%%%%%%%%%%%%%%%%
&
%%%%%%%%%%%%%%%%%%
A_e\boxtimes A_s\dar{\delta\boxtimes \id}&A_s\lar[swap]{\Delta_{e,s}}\rar{\Delta_{s,e}}\dar{\id}&A_s\boxtimes A_e\dar{\id\boxtimes\delta}
\\
 A_r\boxtimes A_{st}\rar{\id_r\boxtimes \Delta_{s,t}}&A_{r}\boxtimes A_s\boxtimes A_t,
%%%%%%%%%%%%%%%%%%%
&
%%%%%%%%%%%%%%%%%%
E\boxtimes A_s\rar{l_{A_s}}&A_s&A_s\boxtimes E\lar[swap]{r_{A_s}}
\end{tikzcd}
\] hold. The composition of two cotensor functors is again a cotensor functor in the sense of \Cref{sec:categ-indep-1:cotensor-functor-composition}. This immediately implies that a cotensor functor $(\mathcal F,\mathcal D,d)$ maps a comonoidal system $(A_s)_{s\in\mathbb S}$ with coproduct morphisms $\Delta_{s,t}$ and counit morphism $\delta$ to a comonoidal systems $(\mathcal F(A_s))_{s\in\mathbb S}$ with coproduct morphisms $\mathcal D_{A_s,A_t}\circ \mathcal F(\Delta_{s,t})$ and counit morphism $d\circ\mathcal F(\delta)$. The analogous statements hold for monoidal systems.

In the following we will only work with comonoidal systems, but since a comonoidal system in $\mathcal C$ is just a monoidal system in $\mathcal C^{\rm op}$, the results have obvious translations for monoidal systems.

\begin{Theorem}\label{sec:comonoidal-systems-1:over-submonoid}
 Let $\mathbb{U}\subset \mathbb S$ be a submonoid. If $(A_s)_{s\in\mathbb S}$ is a comonoidal system with coproduct morphisms $(\Delta_{s,t})_{s,t\in\mathbb S}$ and counit morphism $\delta$, then $(A_s)_{s\in\mathbb{U}}$ is a comonoidal system with coproduct morphisms $(\Delta_{s,t})_{s,t\in\mathbb{U}}$ and counit morphism $\delta$.
\end{Theorem}

\begin{proof} The inclusion $\mathbb{U}\hookrightarrow\mathbb S$ is a monoid homomorphism, hence it is a cotensor functor with respect to the identity natural transformation and identity morphism. The theorem now follows from \Cref{sec:categ-indep-1:cotensor-functor-composition}.
\end{proof}

\subsection{Categorial L\'evy processes}\label{sec:CLP}

\begin{Definition}\label{def:left-divisibility-preorder} Let $\mathbb S$ be a monoid. For $s,t\in\mathbb S$ we write $s\leq t$ if $s$ is a left divisor of $t$, that is if there exists a $r\in\mathbb S$ such that $sr=t$.
\end{Definition}

On every monoid the defined binary relation $\leq$ is a preorder, that is it is reflexive and transitive.

Recall that a monoid $\mathbb S$ is called \emph{cancellative} if $ab=ac$ implies $b=c$ and $ba=ca$ implies $b=c$ for all $a\in\mathbb S$.
In particular, if $s\leq t$ holds in a cancellative monoid $\mathbb S$, then the element $r$ with $sr=t$ is uniquely determined.

\begin{Definition}
 Let $\mathbb S$ be a cancellative monoid. For $s,t\in\mathbb S$ with $s\leq t$, we write $(s\rightarrow t)$ for the (necessarily unique) element with $t=s (s\rightarrow t)$; later we will simply write, in additive notation, $t-s$ instead of $(s\rightarrow t)$.
\end{Definition}

Note that left invertibility, right invertibility and invertibility are all equivalent for elements of a cancellative monoid $\mathbb S$. Indeed, suppose $ab=e$ with $e\in \mathbb S$ the unit element. This implies $baba=bea=bae$. Since $\mathbb S$ is cancellative it follows that $ba=e$ and hence $a=b^{-1}$.

The following definition shall grasp the essence of stationary independent increment processes, cf.\ the remark below.

\begin{Definition}\label{sec:second-induct-limit-2:Def:abstr-LP}
 Let $(A_t)_{t\in\mathbb S}$ be a comonoidal system in $\mathcal C$ over a cancellative monoid~$\mathbb S$. A~\emph{categorial L\'evy process on $(A_t)_{t\in\mathbb S}$} is a collection of morphisms $j_{s,t}\colon A_{(s\rightarrow t)}\to B$ for $s\leq t$ to some common object $B\in\mathcal C$ such that
 \begin{enumerate}\itemsep=0pt
 \item[1)] $j_{t,t}=\Eins_B\circ\delta$,
 \item[2)] $j_{s_1,t_1},\dots, j_{s_n,t_n}$ are independent if $s_1\leq t_1\leq s_2\leq \cdots\leq s_n\leq t_n$,
 \item[3)] $j_{r,s,t}\circ \Delta_{(r\rightarrow s),(s\rightarrow t)}=j_{r,t}$ for some independence morphism $j_{r,s,t}$ of $j_{r,s}$ and $j_{s,t}$.
 \end{enumerate}
\end{Definition}

\begin{Remark}We think of the $j_{s,t}$ as increment of a stationary process over the time interval from $s$ to $t$ and of $A_{r}$ as the common distribution of all $j_{s,t}$ with $(s\rightarrow t)=r$. The required conditions can than be paraphrased as follows. Increments over trivial intervals should vanish; increments over consecutive intervals shall be independent; combining increments over adjacent intervals should yield the increment corresponding to the union of the intervals. In Section~\ref{sec:prob-quant-prob} we will discuss how convolution semigroups of probability measures indeed give rise to (co)monoidal systems in the (opposite) category of probability spaces.
\end{Remark}

In the remainder of this section, we exhibit general conditions on the structure of the index set $\mathbb S$ and the category $\mathcal{C}$ which ensure that comonoidal systems embed into full comonoidal systems (Section~\ref{sec:first-induct-limit}) and that every comonoidal system allows for the construction of a canonical categorial L\'evy process on it (Section~\ref{sec:second-induct-limit}).

\subsection{Unique factorization monoids}\label{sec:uf-monoids}

For an arbitrary set $X$, we denote by $X^*$ the set of all tuples $(x_1,\dots, x_n)$ over $X$ of arbitrary length $n\in\mathbb N_0$. The concatenation of tuples is written in this section as
\begin{align*}
 (x_1,\dots, x_n)\smile (y_1,\dots, y_m):=(x_1,\dots, x_n,y_1,\dots, y_m)
\end{align*}
to clearly distinguish it from other multiplications (later $X$ will typically be a semigroup itself). $X^*$ is a monoid, the \emph{free monoid over $X$}, its neutral element is the empty tuple $()$.

Now suppose that $\mathbb S$ is a monoid with neutral element $e$. A tuple $(s_1,\dots, s_n)\in (\mathbb S\setminus\{e\})^*$ is called a \emph{factorization} of $t\in\mathbb S$ if $t=s_1\cdots s_n$; by convention, for $n=0$ we interpret the empty product on the right hand side as $e$ and, accordingly, the empty tuple is a factorization of $e$. The set of all factorizations of $t\in\mathbb S$ is denoted by $\mathbb{F}_t$. A factorization $\sigma\in\mathbb{F}_t$ is said to be a \emph{refinement} of a factorization $(t_1,\dots, t_n)\in\mathbb{F}_t$ if $\sigma=\tau_1\smile\cdots\smile\tau_n$ for some $\tau_k\in\mathbb{F}_{t_k}$. We write $\sigma\geq\tau$ if $\sigma$ is a refinement of $\tau$. This defines a partial order on $\mathbb{F}_t$.\footnote{\label{fn:J_t-1}Shalit and Skeide have considerations very similar to ours for the sets $\mathbb J_t:=\{(s_1,\dots ,s_n)\in\mathbb S^*\mid t=s_1\dots s_n\}$ whenever $\mathbb S$ is a semigroup (note that they do not exclude $e$ from the possible factors when $\mathbb S$ is a monoid as we do in the definition of $\mathbb F_t$). Because the length of refinements is increasing, it is easy to see that refinement is indeed a partial order, see the proof of \cite[Proposition 4.12]{ShSk20p} for a detailed argument (which applies analogously to $\mathbb F_t$).}

We say that a monoid $\mathbb S$ is \emph{conical} if the only invertible element of $\mathbb S$ is its unit element~$e$. If~$\mathbb S$ is conical, the empty tuple $()$ is the unique factorization of $e$.

\begin{Proposition}\label{sec:cancellative-monoids:prop:unique}
Let $\mathbb S$ be a cancellative, conical monoid. If $\tau_1\smile\cdots\smile\tau_n=\tau'_1\smile\cdots\smile\tau'_n$ with $\tau_k,\tau'_k\in \mathbb{F}_{t_k}$, then $\tau_k=\tau'_k$ for all $k$.
\end{Proposition}

\begin{proof} For $n=0$ or $n=1$ there is nothing to prove. Suppose $\tau_1\smile\cdots\smile\tau_n=\tau'_1\smile\cdots\smile\tau'_n =(s_1,\dots, s_\ell)$ with $\tau_n=(s_k,\dots, s_\ell),\tau'_n=(s_{k'},\dots, s_{\ell})\in\mathbb{F}_{t_n}$. We have $s_k\cdots s_\ell=s_{k'}\cdots s_{\ell}$. Suppose $k<k'$. Since $\mathbb S$ is cancellative, this implies $s_{k}\cdots s_{k'-1}=e$ and thus $s_{k}$ is invertible which contradicts the fact that $\mathbb S$ is conical. So $k\geq k'$. Analogously, we get $k'\geq k$ which shows $k=k'$ and thus $\tau_n=\tau'_n$. Now the proposition follows by induction.
\end{proof}

\begin{Definition}\label{sec:cancellative-monoids:Def:uf-monoid}
 A cancellative monoid $\mathbb S$ is called a \emph{unique factorization monoid} or \emph{uf-monoid} for short if any two factorizations of the same element have a common refinement. A \emph{cuf-monoid} is a conical uf-monoid.
\end{Definition}

Equivalently, a uf-monoid is a cancellative monoid such that $\mathbb{F}_t$ is a directed set with respect to refinement for every $t\in\mathbb S$. The term \emph{unique factorization monoid} or \emph{uf-monoid} has been coined by Johnson \cite{Joh71} who establishes some general theory of uf-monoids, gives different characterizations of the uf-property and presents constructions to find uf-monoids; be aware that the term can be slightly misleading, because typically there are no (finite) unique factorizations of elements in a uf-monoid. We only deal with cuf-monoids, because we will need \Cref{sec:cancellative-monoids:prop:unique}. The examples we will use later on are only $\mathbb N_0$, $\mathbb Q_+$ and $\mathbb R_+$ (with addition), but it seems that cuf-monoids provide the most general setting in which we can study the inductive limit constructions of the following sections.

\subsection{First inductive limit: the generated full comonoidal system}\label{sec:first-induct-limit}

 Let $((A_s)_{s\in\mathbb S},(\Delta_{s,t})_{s,t\in\mathbb S},\delta)$ be a comonoidal system over a cancellative and conical monoid $\mathbb S$ in a tensor category $(\mathcal C,\boxtimes)$ with unit object $E$. For a tuple $\sigma=(s_1,\dots,s_n)\in\mathbb{F}_t$ put $A_{\sigma}:=A_{s_1}\boxtimes\cdots\boxtimes A_{s_n}$ for $n\geq 1$ and $A_{()}:=E$. Define $\Delta_{\sigma}\colon A_t\to A_\sigma$ recursively by
 \begin{gather*}
 \Delta_{()}:=\delta,\\
 \Delta_{(t)}:=\id_{t},\\
 \Delta_{(s_1,\dots, s_{n+1})}:=(\Delta_{(s_1,\dots,s_n)}\boxtimes\id_{s_{n+1}})\circ\Delta_{s_1\cdots s_n,s_{n+1}}\quad\text{for $n\geq1$;}
 \end{gather*}
 note that, for $n=1$, $\Delta_{(s_1,s_2)}=(\id_{s_1}\boxtimes\id_{s_2})\circ\Delta_{s_1,s_2}=\Delta_{s_1,s_2}$ gives back the original binary coproduct morphisms and that the third equation actually also holds for $n=0$ by the counit property, $(\delta\boxtimes \id_{s})\circ\Delta_{e,s}=\id_s$ (suppressing the unit constraint).
Let $\tau=(t_1,\dots,t_n)\in\mathbb{F}_t$ and $\sigma\geq\tau$. Since $\mathbb S$ is cancellative and conical, we can use \Cref{sec:cancellative-monoids:prop:unique} to write $\sigma=\tau_1\smile\cdots\smile\tau_n$ for uniquely determined $\tau_k\in\mathbb{F}_{t_k}$. With this notation we put
\begin{align*}
 \Delta_{\sigma}^\tau:=\Delta_{\tau_1}\boxtimes \cdots\boxtimes \Delta_{\tau_n}\colon \ A_\tau\to A_\sigma
\end{align*}
for all $\tau\leq\sigma\in\mathbb{F}_t$.

\begin{Lemma}\label{CS:Lem:ind-syst1} It holds that $\Delta_\rho^\sigma\circ\Delta_{\sigma}^\tau=\Delta_\rho^\tau$
 for all $\tau\leq\sigma\leq\rho\in\mathbb{F}_t$ for every cancellative and conical monoid.
\end{Lemma}

\begin{proof} The proof of Bhat and Mukherjee for the case $\mathbb S=\mathbb R_+$ \cite[Lemma~4]{BhMu10} works without a~change for general cancellative and conical monoids.
\end{proof}

\begin{Corollary} Let $\mathbb S$ be a cuf-monoid and $\big((A_s)_{s\in\mathbb S},(\Delta_{s,t})_{s,t\in\mathbb S},\delta\big)$ a~comonoidal system in a~tensor category $\mathcal C$. Then for every $t\in\mathbb S$, $\bigl((A_\tau)_{\tau\in\mathbb{F}_t},(\Delta_{\sigma}^{\tau})_{\sigma\geq\tau\in\mathbb{F}_t}\bigr)$ is an inductive system.
\end{Corollary}

\begin{proof} By \Cref{sec:cancellative-monoids:Def:uf-monoid} of a uf-monoid, $\mathbb{F}_t$ is directed. The first condition of an inductive system, $\Delta_{\tau}^{\tau}=\id_{\tau}$, is obvious. The second condition, $\Delta_\rho^\sigma\circ\Delta_{\sigma}^\tau=\Delta_\rho^\tau $ for all $\tau\leq\sigma\leq\rho\in\mathbb{F}_t$, is the statement of \Cref{CS:Lem:ind-syst1}.
\end{proof}

Suppose that the inductive systems $(A_\tau)_{\tau\in\mathbb{F}_t}$ have inductive limits $\mathscr A_t$ with morphisms $D^\tau\colon A_\tau\allowbreak \to\mathscr A_t$. For $\tau\in\mathbb{F}_t$ denote by $\mathbb{F}_\tau$ the set of all refinements of $\tau$. Then $\mathbb{F}_\tau$ is a cofinal subset of~$\mathbb{F}_t$, see \Cref{sec:inductive-limits-1:ex:cofinal}. Denote by $\mathscr A_\tau$ the inductive limit. Then there is a canonical isomorphism $\mathscr A_t\cong\mathscr A_\tau$ because of \Cref{sec:inductive-limits-2:cor:cofinal}.

\begin{Lemma} The diagram
\[
 \begin{tikzcd}[column sep=huge]
 A_{\sigma}\boxtimes A_{\tau}\rar{D^{\sigma}\boxtimes D^{\tau}} \dar[swap]{\Delta_{\sigma'}^\sigma\boxtimes\Delta_{\tau'}^\tau}&\mathscr A_s\boxtimes\mathscr A_t\\
 A_{\sigma'}\boxtimes A_{\tau'}\urar[swap]{D^{\sigma'}\boxtimes D^{\tau'}}&
 \end{tikzcd}
\]
 commutes for all $\sigma'\geq\sigma\in\mathbb{F}_s$, $\tau'\geq\tau\in\mathbb{F}_t$.
\end{Lemma}

\begin{proof}
 By functoriality of $\boxtimes$, we have
 \[
 \big(D^{\sigma'}\boxtimes D^{\tau'}\big)\circ\big(\Delta_{\sigma'}^\sigma\boxtimes\Delta_{\tau'}^\tau\big)
 =\big(D^{\sigma'}\circ\Delta_{\sigma'}^{\sigma}\big)\boxtimes\big(D^{\tau'}\circ\Delta_{\tau'}^\tau\big)
 =D^{\sigma}\boxtimes D^{\tau}.\tag*{\qed}
 \]\renewcommand{\qed}{}
\end{proof}

So, by the universal property of the inductive limit, there are unique morphisms $\widetilde{\Delta}_{s,t}\colon \mathscr A_{st}\to\mathscr A_s\boxtimes\mathscr A_t$ such that
 \[
 \begin{tikzcd}[column sep=huge]
 A_{\sigma}\boxtimes A_{\tau}\rar{D^{\sigma}\boxtimes D^{\tau}} \dar{D^{\sigma\smile\tau}}&\mathscr A_s\boxtimes\mathscr A_t\\
 \mathscr A_{st}\cong\mathscr A_{(s,t)}\urar[swap]{\widetilde{\Delta}_{s,t}}&
 \end{tikzcd}
 \]
 commutes for every $\sigma\in\mathbb{F}_s$, $\tau\in\mathbb{F}_t$.

\begin{Theorem}\label{sec:first-induct-limit-1:theo:comon-syst}
 The $\mathscr A_t$ form a comonoidal system with respect to the coproduct morphisms~$\widetilde{\Delta}_{s,t}$ and the counit morphism~$\id_E$.
\end{Theorem}

\begin{proof}
 Note that we identify the (trivial) inductive limit $\mathscr A_e$ with $E$, which is feasible because $\mathbb{F}_e=\{()\}$ and $A_{()}=E$. The counit property is then trivially fulfilled. In the diagram
\[
\begin{tikzcd}[column sep=huge]
 &\mathscr A_{rst}\dlar[swap]{\widetilde\Delta}\drar{\widetilde\Delta}&\\
 \mathscr A_{r}\boxtimes\mathscr A_{st}\drar[swap]{\id\boxtimes\widetilde\Delta}&A_{\rho}\boxtimes A_{\sigma}\boxtimes A_{\tau}\uar{D^{\rho\smile\sigma\smile\tau}}\rar{D^{\rho\smile\sigma}\boxtimes D^{\tau}}\lar[swap]{D^{\rho}\boxtimes D^{\sigma\smile\tau}}\dar{D^\rho\boxtimes D^\sigma\boxtimes D^\tau}&\mathscr A_{rs}\boxtimes\mathscr A_t\dlar{\widetilde{\Delta}\boxtimes\id}\\
 &\mathscr A_r\boxtimes\mathscr A_s\boxtimes\mathscr A_t&
\end{tikzcd}
\]
the four corners commute by the definition of $\widetilde{\Delta}$. From $(r,s,t)\in\mathbb{F}_{rst}$ we get a canonical isomorphism $\mathscr A_{rst}\cong\mathscr A_{(r,s,t)}$ which we use to identify the two. Coassociativity now follows from the universal property of the inductive limit.
\end{proof}

Define $D_t\colon A_t\to\mathscr A_t$ by $D_t:=D^{(t)}$ for $t\neq e$ and $D_e:=\delta$.

\begin{Theorem}\label{sec:first-induct-limit-1:theo:emb-into-full}
 The morphisms $(D_t)_{t\in\mathbb S}$ form a \emph{morphism of comonoidal systems}, that is $\widetilde\Delta_{s,t}\circ D_{st}=(D_{s}\boxtimes D_{t})\circ \Delta_{s,t}$ and $\id_E\circ D_{e}=\delta$.
\end{Theorem}

\begin{proof}
 The counit is respected by definition of $D_e$. In the diagram
\[
\begin{tikzcd}
 A_{st}\rar{\Delta_{s,t}}\dar{D^{(st)}}&A_s\boxtimes A_t\dlar{D^{(s,t)}}\dar{D^{(s)}\boxtimes D^{(t)}}\\
 \mathscr A_{st}\cong\mathscr A_{(s,t)}\rar[swap]{\widetilde\Delta_{s,t}}&\mathscr A_s\boxtimes\mathscr A_t
\end{tikzcd}
\]
the lower right commutes by definition of $\widetilde\Delta$ and the upper left because $\mathscr A_{st}$ is the inductive limit. So the outside square commutes, which finishes the proof.
\end{proof}

Let $\mathcal F\colon \mathcal C\to\mathcal D$ be a functor. Then any inductive system $\bigl((A_\alpha)_\alpha,(f^\alpha_\beta)_{\alpha\leq\beta}\bigr)$ in $\mathcal C$ is mapped to an inductive system $\bigl((\mathcal F(A_\alpha))_{\alpha},(\mathcal F(f^\alpha_\beta))_{\alpha\leq\beta}\bigr)$ in $\mathcal D$.
We say that $\mathcal F$ \emph{preserves inductive limits} if for every inductive limit $\bigl(\mathscr A,(f_{\alpha})_\alpha\bigr)$ of an inductive system $\bigl((A_\alpha)_\alpha,(f^\alpha_\beta)_{\alpha\leq\beta}\bigr)$ it holds that $\bigl(\mathcal F(\mathscr A),(\mathcal F(f_{\alpha}))_{\alpha}\bigr)$ is an inductive limit of the inductive system $\bigl((\mathcal F(A_\alpha))_{\alpha},(\mathcal F(f^\alpha_\beta))_{\alpha\leq\beta}\bigr)$.

\begin{Theorem}\label{sec:first-induct-limit-1:theo:full-system}
 If the tensor product preserves inductive limits, the morphisms $\widetilde{\Delta}_{s,t}$ are all isomorphisms. In other words $(\mathscr A_t)_{t\in\mathbb S}$ is a full comonoidal system.
\end{Theorem}

\begin{proof}
 The tensor product is a bifunctor $\boxtimes\colon \mathcal C\times\mathcal C\to\mathcal C$. Inductive systems in $\mathcal C\times\mathcal C$ are in bijection with pairs of inductive systems in $\mathcal C$ and an inductive limit in $\mathcal C\times\mathcal C$ is a pair of inductive limits for the inductive systems in $\mathcal C$. If $\boxtimes$ preserves inductive limits, $\mathscr A_s\boxtimes\mathscr A_t$ is an inductive limit of the inductive system formed by $(A_\sigma\boxtimes A_\tau)_{\sigma\in\mathbb{F}_s,\tau\in\mathbb{F}_t}$ with respect to the maps $D^{\sigma}\boxtimes D^{\tau}$. Since $\widetilde{\Delta}_{s,t}$ makes the diagram
\[
 \begin{tikzcd}[column sep=huge]
 A_{\sigma}\boxtimes A_{\tau}\rar{D^{\sigma}\boxtimes D^{\tau}} \dar[swap]{D^{\sigma\smile\tau}}&\mathscr A_s\boxtimes\mathscr A_t\\
 \mathscr A_{st}\cong\mathscr A_{(s,t)}\urar[swap]{\widetilde{\Delta}_{s,t}}&
 \end{tikzcd}
 \]
commute, it is the canonical isomorphism between the two inductive limits.
\end{proof}

All tensor categories we are interested in have tensor products which do preserve inductive limits. % The following example illustrates that this is not true in general.

\subsection{Ore monoids}\label{sec:ore-monoids}

\begin{Definition} A cancellative monoid $\mathbb S$ is called \emph{Ore monoid} if for all $s,t\in\mathbb S$ there exists an $r\in\mathbb S$ with $s\leq r $ and $t\leq r$.

In other words: A cancellative monoid is an Ore monoid if and only if $(\mathbb S,{\leq})$ is directed.
\end{Definition}

Be aware of our choice to use left divisibility in the definition while some authors use the same term for the symmetric notion based on right divisibility instead.

Product systems of C*-correspondences (Hilbert bimodules) over Ore monoids have been studied independently by Albandik and Meyer~\cite{AlMe15} and Kwa{\'s}niewski and Szyma{\'n}ski~\cite{KwSz16}. Albandik and Meyer's definitions are slightly more general, because they allow also noncancellative monoids. This works as well, but one has to replace inductive limit by filtered colimits which many readers are probably less farmiliar with. The recent work of Shalit and Skeide~\cite{ShSk20p} on subproduct systems also emphasizes the importance of the Ore condition in this context.

We could not find the following theorem in the literature, and it seems that it is at least not well known in the product system community. The backward direction has recently been proved independently by Shalit and Skeide as \cite[Theorem~4.13]{ShSk20p} (actually they prove a slightly more general statement because they do not restrict to conical monoids).

\begin{Theorem}\label{thm:totally_ordered-iff-oreuf} Let $\mathbb S$ be a cancellative monoid. Then the following are equivalent:
 \begin{itemize}\itemsep=0pt
 \item $\mathbb S$ is a cuf-monoid and an Ore monoid.
 \item $\mathbb S$ is totally ordered with respect to $\leq$.
 \end{itemize}
\end{Theorem}

\begin{proof}
 Suppose that $\mathbb S$ is a cuf-monoid and an Ore monoid. Because $\mathbb S$ is an Ore monoid, $(\mathbb S,\leq)$ is a directed set. Let $s,t\in\mathbb S$. By the Ore property there exist $r,p,q\in\mathbb S$ with $r=sp=tq$. We want to show that $s\in t\mathbb S$ or $t\in s\mathbb S$. If one of the four elements $s$, $t$, $p$, $q$ is $e$ this is obvious. If none of them is invertible, then $(s,p),(t,q)\in \mathbb{F}_r$, so they have a common refinement $(r_1,\dots, r_n)\in\mathbb{F}_r$ because $\mathbb S$ is a uf-monoid. This means $s=r_1\cdots r_i$, $q=r_{i+1}\cdots r_n$ and $t=r_1\cdots r_j$, $q=r_{j+1}\cdots r_n$ for some $i,j\in\{1,\dots,n\}$. Now $i\leq j$ implies $t=s r_{i+1}\cdots r_j\in s\mathbb S$ and $j\leq i$ implies $s=t r_{j+1}\cdots r_i\in t\mathbb S$, so we are done.

Now suppose that $\mathbb S$ is totally ordered. Clearly, $\mathbb S$ is an Ore monoid. Let $\varepsilon\in U(\mathbb S)$. Then $\varepsilon\leq e$ and $e\leq \varepsilon$ imply $\varepsilon=e$, so there are no nontrivial invertible elements. We prove that for all $r\in\mathbb S$ and all $(t_1,\dots,t_n),(s_1,\dots, s_m)\in\mathbb{F}_r$ there is a common refinement by induction on $n$. For $n=0$ and $n=1$ this is obvious. Now suppose the statement holds for $n-1\in\mathbb N$ and let $(t_1,\dots,t_n),(s_1,\dots, s_m)\in\mathbb{F}_r$. Because $t_1\cdots t_n= r$ we know that $t_1\cdots t_{n-1} < r$. Because $\leq$ is a~total order, we are in one of the following situations:

\emph{Case 1:} $t_1\cdots t_{n-1}= s_1\cdots s_k$ for some $k\in\{1,\dots,m-1\}$. By the induction hypthesis, there is a common refinement $(r_1,\dots, r_\ell)$ of $(t_1,\dots, t_{n-1})$ and $(s_1,\dots ,s_k)$. It is easy to check that $(r_1,\dots,r_\ell,s_{k+1},\dots,s_m)$ is a common refinement of $(t_1,\dots, t_n)$ and $(s_1,\dots, s_m)$.

\emph{Case 2:} $s_1\cdots s_{k}<t_1\cdots t_{n-1}<s_1\cdots s_{k+1}$ for some $k\in\{1,\dots,m-1\}$. Then there are $p$,~$q$ such that $s_1\cdots s_k p=t_1\cdots t_{n-1}$ and $t_1\cdots t_{n-1}q=s_1\cdots s_{k+1}$. By the induction hypothesis, there exists a comon refinement $(r_1,\dots,r_\ell)$ of $(s_1,\dots ,s_k ,p)$ and $(t_1,\dots, t_{n-1})$. Now it follows that $(r_1,\dots,r_\ell,q,s_{k+2},\dots,s_m)$ is a common refinement of $(t_1,\dots,t_n)$ and $(s_1,\dots, s_m)$.
\end{proof}

\begin{Remark}\label{rem:ShSk-question} A natural question is whether a cuf monoid is automatically totally ordered, cf.\ the question Shalit and Skeide ask in \cite[last bullet at the end of Section~4]{ShSk20p}.\footnote{Their question is almost that, but the setup is slightly different. They ask whether a not necessarily conical monoid in which the sets $\mathbb J_t:=\{(s_1,\dots ,s_n)\in\mathbb S^*\mid t=s_1\dots s_n\}$ (cf.\ footnote~\ref{fn:J_t-1}) are all directed is automatically cancellative and \emph{totally directed}, which is a generalization of totally ordered, not demanding divisibility to be an antisymmetric relation. It is easy to see that directedness of $\mathbb J_t$ and $\mathbb F_t$ are equivalent for every monoid $\mathbb S$. Assuming also that $\mathbb S$ is conical, their question translates to the one we posed.} Johnson~\cite{Joh71} proves that free products of uf-monoids are again uf. In particular, every free monoid on more than one generator gives an example of a cuf monoid which is (quite obviously) not Ore and, therefore, is not totally directed. This resolves the question in the negative.
\end{Remark}

\subsection{Second inductive limit: L\'evy processes}\label{sec:second-induct-limit}

Let $\mathscr A_t$ be a full comonoidal system over an Ore monoid $\mathbb S$ with coproduct isomorphisms $\widetilde{\Delta}_{s,t}$ in a tensor category with initial unit object. Without loss of generality assume that $\mathscr A_e=E$ and that the counit morphism is $\delta=\id_E$.

For $s\leq t$, $t=sp$, define $i^s_t\colon \mathscr A_s\to\mathscr A_t$ as the composition
\[
\begin{tikzcd}
 \mathscr A_s\rar{\iota^1}&\mathscr A_s\boxtimes\mathscr A_{p}\rar{\widetilde{\Delta}^{-1}}&\mathscr A_t.
\end{tikzcd}
\]

\begin{Theorem} $\big((\mathscr A_t)_{t\in\mathbb S},(i^{s}_t)_{s\leq t}\big)$ is an inductive system.
\end{Theorem}

\begin{proof} Let $r\leq s\leq t$, $s=rp$, $t=sq$. In the diagram
\[
\begin{tikzcd}[column sep=huge]
 \mathscr A_r\dar{\id}\rar{\iota^1}&\mathscr A_{r}\boxtimes\mathscr A_{p}\dar[swap]{\id\boxtimes\iota^1}\dar{\iota^1}\rar{\widetilde\Delta^{-1}}& \mathscr A_{s}\dar{\iota^1} \\
 \mathscr A_r \dar{\id}\rar{\iota^1}&\mathscr A_r\boxtimes\mathscr A_{p}\boxtimes\mathscr A_{q} \dar{\id\boxtimes\widetilde\Delta^{-1}} \rar{\widetilde\Delta^{-1}\boxtimes\id} &\mathscr A_{s}\boxtimes\mathscr A_{q}\dar{\widetilde\Delta^{-1}} \\
 \mathscr A_r\rar{\iota^1}&\mathscr A_r\boxtimes\mathscr A_{pq} \rar{\widetilde\Delta^{-1}} & \mathscr A_t
\end{tikzcd}
\]
the lower right corner commutes by coassociativity of $\widetilde\Delta$ and the other three corners commute by the naturality of $\iota^1$. We suppressed the associativity constraint and identified $\mathscr A_r\boxtimes(\mathscr A_{p}\boxtimes\mathscr A_{q})$ with $(\mathscr A_r\boxtimes\mathscr A_{p})\boxtimes\mathscr A_{q}$, which leads to the two interpretations $\id_{\mathscr A_r}\boxtimes\iota^1_{\mathscr A_{p},\mathscr A_{q}}$ and $\iota^1_{\mathscr A_r\boxtimes\mathscr A_{p},\mathscr A_{q}}$ of the arrow $\mathscr A_{r}\boxtimes\mathscr A_{p}\to\mathscr A_r\boxtimes\mathscr A_{p}\boxtimes\mathscr A_{q}$. Of course, this would not work without coherence as discussed in \Cref{obs:coherence}.
\end{proof}

For the rest of this section, we fix a totally ordered monoid $\mathbb S$ and an inductively complete tensor category $\mathcal C$ with compatible inclusions $\iota^1$, $\iota^2$, whose tensor product $\boxtimes$ perserves inductive limits.
\begin{Remark}
 As the typical examples of totally ordered monoids are submonoids of $\mathbb R_+$, we will use additive notation in the following. For $s\leq t$, we denote the unique element $r$ with $s+r=t$ by $t-s$.
\end{Remark}

Given only the comonoidal system $\bigl((A_t)_{t\in\mathbb S},(\Delta_{s,t})_{s,t\in\mathbb S}\bigr)$ one can construct a canonical categorial L\'evy process. Since $\mathbb S$ is a uf-monoid and $\boxtimes$ preserves inductive limits, $(A_t)_{t\in\mathbb S}$ generates a full comonoidal system $\bigl((\mathscr A_t)_{t\in\mathbb S},\widetilde\Delta\bigr)$ by \Cref{sec:first-induct-limit-1:theo:full-system}. Denote by $D^{t}\colon A_t\to\mathscr A_t$ the canonical morphisms. Let $\bigl(\mathscr A,(i^t\colon \mathscr A_t\to\mathscr A)_{t\in\mathbb S}\bigr)$ be the inductive limit of $(\mathscr A_t)_{t\in\mathbb S}$. Define $j_{s,t}\colon A_{t-s}\to\mathscr A$ as the composition
\[
\begin{tikzcd} A_{t-s}\rar{D^{t-s}}&\mathscr A_{t-s}\rar{\iota^2}&\mathscr A_s\boxtimes\mathscr A_{t-s}\rar{\widetilde\Delta^{-1}}&\mathscr A_t\rar{i^t}&\mathscr A.
\end{tikzcd}
\]

\begin{Theorem}\label{sec:second-induct-limit-1:theo-LP}
 The $j_{s,t}$ form a categorial L\'evy process.
\end{Theorem}

\begin{proof}
 We construct the independence morphism $j_{r,s,t}$ for $j_{r,s}$ and $j_{s,t}$ and show that $j_{r,s,t}\circ\Delta_{s-r,t-s}=j_{r,t}$. Define $j_{r,s,t}$ as the composition
\[
\begin{tikzcd}
 A_{s-r}\boxtimes A_{t-s}\rar{D^{s-r}\boxtimes D^{t-s}}&\mathscr A_{s-r}\boxtimes\mathscr A_{t-s}\rar{\iota^2}&\mathscr A_r\boxtimes\mathscr A_{s-r}\boxtimes\mathscr A_{t-s}\rar{\widetilde\Delta^{-1}}&\mathscr A_t\rar{i^t}&\mathscr A.
\end{tikzcd}
\]
Now, the diagram
\[
\begin{tikzcd}
 A_{t-s}\rar{D^{t-s}}\dar{\iota^2}&\mathscr A_{t-s}\rar{\iota^2}\dar{\iota^2}&\mathscr A_s\boxtimes\mathscr A_{t-s}\rar{\widetilde\Delta^{-1}}\dar{\widetilde\Delta}&\mathscr A_t\drar{i^t}\dar{\id}&\\
 A_{s-r}\boxtimes A_{t-s}\rar{D^{s-r}\boxtimes D^{t-s}}&\mathscr A_{s-r}\boxtimes\mathscr A_{t-s}\rar{\iota^2}&\mathscr A_r\boxtimes\mathscr A_{s-r}\boxtimes\mathscr A_{t-s}\rar{\widetilde\Delta^{-1}}&\mathscr A_t\rar{i^t}&\mathscr A\\
A_{s -r}\rar{D^{s -r}}\uar{\iota^1}&\mathscr A_{s -r}\rar{\iota^2}\uar{\iota^1}&\mathscr A_r\boxtimes\mathscr A_{s -r}\rar{\widetilde\Delta^{-1}}\uar{\iota^1}&\mathscr A_s \uar{i^{s}_t}\urar{i^s }&
\end{tikzcd}
\]
commutes: The leftmost squares commute due to naturality of $\iota^1$ and $\iota^2$. The next squares commute by \Cref{obs:coherence} and naturality of $\iota^2$. The upper right square commutes by coassociativity of $\widetilde\Delta$. The triangles commute by definition of the inductive limit. It remains to show commutativity of the lower right square. In a bit more detail, this is
\[
\begin{tikzcd}[column sep=large] \mathscr A_r\boxtimes\mathscr A_{s-r}\boxtimes\mathscr A_{t-s}\rar{\widetilde\Delta^{-1}\boxtimes\id}&\mathscr A_s\boxtimes\mathscr A_{t-s}\rar{\widetilde\Delta}&\mathscr A_t,\\
\mathscr A_r\boxtimes\mathscr A_{s -r}\rar{\widetilde\Delta^{-1}}\uar{\iota^1}&\mathscr A_s\uar{\iota^1}\urar{i^s_{t}}
\end{tikzcd}
\]
which commutes by naturality of $\iota^1$ and the definition of $i^{s}_{t}$. This shows that $j_{r,s,t}$ is an independence morphism. Next we consider
\[
\begin{tikzcd}
 A_{t-r}\rar{D^{t-r}}\dar{\Delta_{s-r,t-s}}&\mathscr A_{t-r}\rar{\iota^2}\dar{\widetilde\Delta}&\mathscr A_r\boxtimes\mathscr A_{t-r}\rar{\widetilde\Delta^{-1}}\dar{\id\boxtimes\widetilde\Delta}&\mathscr A_t\drar{i^t}\dar\id&\\
 A_{s-r}\boxtimes A_{t-s}\rar{D^{s-r}\boxtimes D^{t-s}}&\mathscr A_{s-r}\boxtimes\mathscr A_{t-s}\rar{\iota^2}&\mathscr A_r\boxtimes\mathscr A_{s-r}\boxtimes\mathscr A_{t-s}\rar{\widetilde\Delta^{-1}}&\mathscr A_t\rar{i^t}&\mathscr A
\end{tikzcd}
\]
in which the first square commutes because the $D^{t}$ form a morphism of comonoidal systems by \Cref{sec:first-induct-limit-1:theo:emb-into-full}, the second square commutes due to naturality of $\iota^2$, and the last square and the triangle commute trivially. So the outside commutes, thus establishing $j_{r,s,t}\circ\Delta_{s-r,t-s}=j_{r,t}$.

The general construction of an independence morphism for $j_{t_1,t_2},\dots, j_{t_{n},t_{n+1}}$ works similar to that for $j_{r,s}$ and $j_{s,t}$.
\end{proof}

There is also a direct way from the comonoidal system $(A_t)_{t\in\mathbb S}$ to the L\'evy process. Put $\mathbb{F}:=\{\sigma=(s_1,\dots, s_n)\mid n\in\mathbb N_0, s_k\in\mathbb S\setminus\{e\}\}=\bigcup_{s\in\mathbb S}\mathbb{F}_s$ and define $\sigma\geq\tau=(t_1,\dots, t_n)$ if there exist $\tau_1,\dots, \tau_n,\tau_{n+1}$ with $\sigma=\tau_1\smile\cdots\smile\tau_n\smile\tau_{n+1}$, $\tau_k\in \mathbb{F}_{t_k}$ for $k\in\{1,\dots, n\}$ and $\tau_{n+1}\in\mathbb{F}$. One shows that $\mathbb{F}$ is directed analogously to $\mathbb{F}_t$. Then define an inductive system $(A_{\sigma})_{\sigma\in\mathbb{F}}$ with respect to the morphisms $i^{\tau}_\sigma\colon A_\tau\to A_\sigma$ defined as the composition
\[
\begin{tikzcd} A_{\tau}\ar{rr}{\Delta^{\tau}_{\tau_1\smile\cdots\smile\tau_n}}&&A_{\tau_1}\boxtimes\cdots\boxtimes A_{\tau_n}\rar{\iota^1}&A_{\tau_1}\boxtimes\cdots\boxtimes A_{\tau_n}\boxtimes A_{\tau_{n+1}}=A_\sigma
\end{tikzcd}
\]
for $\sigma=\tau_1\smile\cdots\smile\tau_n\smile\tau_{n+1}\geq\tau$.

\begin{Theorem}\label{theo:can-iso}
 The inductive limits of $(A_{\sigma})_{\sigma\in\mathbb{F}}$ and $(\mathscr A_t)_{t\in\mathbb S}$ are canonically isomorphic.
\end{Theorem}

\begin{proof}
 This is exactly the situation of \Cref{prop:ind-lim-agree}.
\end{proof}

\subsection{Summary of the necessary assumptions}
\label{sec:assumptions}

We end this section with a short summary what conditions on $\mathbb S$, $(A_s)_{s\in\mathbb S}$ and $\mathcal C$ are needed for each of the presented steps.
 \begin{enumerate}\itemsep=0pt
 \item For \Cref{sec:second-induct-limit-2:Def:abstr-LP} (categorial L\'evy process), we assume $\mathbb S$ to be cancellative and $\mathcal C$ to have an initial unit object.
 \item For the first inductive limit as described in Section~\ref{sec:first-induct-limit} (full comonoidal system from comonoidal system), we assume that $\mathbb S$ is a cuf-monoid, that the inductive limits $\mathscr A_t$ exist in $\mathcal C$, and that $\boxtimes$ preserves inductive limits.
 \item For the second inductive limit as described in Section~\ref{sec:second-induct-limit} (categorial L\'evy process from full comonoidal system), we need $\mathbb S$ to be an Ore monoid and the inductive limit~$\mathscr A$ to exist in~$\mathcal C$.
 \item For the combination of both inductive limits (categorial L\'evy process from comonoidal system), we need $\mathbb S$ to be both, Ore and cuf; this is equivalent to being totally ordered with respect to the divisibility relation~$\leq$. Of course, the relevant inductive limits have to exist in~$\mathcal C$ and $\boxtimes$ must preserve inductive limits.
 \end{enumerate}

\section{Examples}\label{sec:examples-1}

We explore what categorial independence and comonoidal systems are in several concrete tensor categories. First we shall see that categorial independence encompasses notions such as linear independence and orthogonality. Then we shall see that the quantum probabilistic notions of independence induced by universal products (e.g., Voiculescu's freeness) and newer generalizations are covered and yield quantum L\'evy processes.

A general pattern for finding interesting independences turns out to be the following. Start with a category $C$ with binary coproducts $A\sqcup B$ and an initial object $E$. As discussed before, independence is trivial in this case because $j_1\sqcup\dots\sqcup j_n\colon B_1\sqcup\dots \sqcup B_n\to A$ is always an independence morphism for $j_k\colon B_k\to A$, $k=1,\dots,n$. However, if we move to a subcategory $C'$ by restricting the class of morphisms (making sure that the coproduct of $C$ restricts to a tensor product on $C'$ and that the initial object of $C$ is still initial in $C'$), independence of $j_k\colon B_k\to A$ in $C'$ is equivalent to $\bigsqcup_{k} j_k\colon \bigsqcup_k B_k\to A$ being a morphism in $C'$.

\subsection{Independence in nonprobabilistic categories}\label{sec:nonpr-categ}

{\bf Pairwise disjointness.}
The category $\Set$ has the empty set $\varnothing$ as initial object and disjoint union $\dot\cup$ as coproduct. By $\Set^{inj}$, we denote the category of sets with injective maps, which is a tensor category with respect to disjoint union with initial unit object $\varnothing$. Injective maps $j_k\colon B_k\to A$ are independent if and only if the canonical map $B_1\dot\cup\cdots\dot\cup B_n\to A$ is injective. This is obviously equivalent to pairwise disjointness of the subsets $j_k(B_k)\subset A$. In this sense, pairwise disjointness of subsets is a special case of categorial independence.

{\bf Linear independence.} Consider the category $\Vect^{inj}$ of vector spaces with injective linear maps. The direct sum is the coproduct in $\Vect$ (vector spaces with arbitrary linear maps), and it turns $\Vect^{inj}$ into a tensor category with initial unit object $\{0\}$. Injective linear maps $f_i\colon V_i\to W$, $i=1,\dots,n$ are independent if and only if they have linearly independent ranges; indeed, the only choice for the independence morphism is the linear map $h:=f_1+\dots + f_n\colon V_1\oplus\dots \oplus V_n\to W$. If $h$ is injective, then $f_1(v_1)+\dots + f_n(v_n)=0$ implies $f_1(v_1)=\dots =f_n(v_n)=0$, so the ranges are linearly independent. On the other hand, if the ranges are linearly independent and $h(v_1\oplus\dots \oplus v_n)=0$, we can conclude that $f_i(v_i)=0$ for $i\in\{1,\dots,n\}$. Since the $f_i$ are injections, it follows that $v_i=0$ for all $i$, so $h$ is injective.

{\bf Algebraic independence.}
Let $L\supset K$ be a field extension. Elements $a_1,\dots, a_n\in L$ are algebraically independent if the canonical map from the polynomial algebra $K[x_1,\dots,x_n]$ to $K[a_1,\dots,a_n]\subset L$ is an isomorphism. The category of commutative unital $K$-algebras with unital algebra homomorphisms has the tensor product as coproduct and the one-dimensional algebra $K$ as unit object. When we restrict, as in the previous examples, to injective homomorphisms, we get another tensor category with initial unit object. The inclusion maps $K[a_i]\hookrightarrow L$ are categorially independent if and only if the canonical $K$-algebra homomorphism $K[a_1]\otimes \dots\otimes K[a_n]\to L$ is injective. Injectivity however means that it is an isomorphism onto its image $K[a_1\dots a_n]\subset L$. Assume that $a_1,\dots, a_n$ are transcendent over~$K$, so that $K[a_i]\cong K[x_i]$. Than categorial independence of the embeddings is equivalent to $K[x_1,\dots,x_n]\cong K[x_1]\otimes\dots K[x_n] \cong K[a_1]\otimes \dots\otimes K[a_n]\cong K[a_1,\dots, a_n]$, hence to algebraic independence of the $a_i$.

{\bf Orthogonality.}
Similar to the previous examples, the category $(\Hilb^{\rm isom},\oplus)$ of Hilbert spaces with isometries is a tensor category with initial unit object $\{0\}$. Isometries $v_i\colon H_i\to G$ are independent if and only if they have orthogonal ranges. Indeed, the only choice for the independence morphism is the linear map $h=v_1+\dots +v_n$. This is an isometry if and only it is a unitary onto its range $v_1(H_1) + \dots + v_n(H_n)$, and this is equivalent to $v_1(H_1), \dots v_n(H_n)$ being pairwise orthogonal.

\subsection{(Co)monoidal systems in nonprobabilistic categories}\label{sec:examples-monoidal-systems-nonprob}

{\bf Graded algebras.}
Monoids in $(\Vect,\otimes)$ are unital algebras, comonoids are coalgebras. Let $(A_t)_{t\in\mathbb S}$ be a monoidal system in $(\Vect,\otimes)$. Then $A:=\bigoplus_{t\in\mathbb S} A_t$ is an $\mathbb S$-graded algebra with respect to the multiplication given by
\begin{align*}
 ab:=\mu_{s,t}(a\otimes b)
\end{align*}
for elements $a\in A_s$, $b\in A_t$ (which has also been observed in \cite[Example~3.18]{AgMa10}). If we consider $(\Vect^{\rm surj},\otimes)$, a monoidal system over $\mathbb N_0$ yields a \emph{standard graded algebra}, i.e., an $\mathbb N_0$-graded algebra $A= \bigoplus_{n\in\mathbb N_0} A_n$ with $A_0=\mathbb C\Eins$ and $A_mA_n=A_{m+n}$. Indeed, the two conditions are exactly the surjectivity of the unit morphism $\Eins\colon \mathbb C\to A_0$ and the product morphisms $\mu_{m,n}\colon A_m\otimes A_n\to A_{m+n}$.

{\bf Subproduct systems.}
Comonoidal systems in $(\Hilb^{\rm isom},\otimes)$ with $H_e=\mathbb C$ and $\delta=\id_\mathbb C$ are called \emph{subproduct systems}. Full subproduct systems are called \emph{product systems}. In this context, the inductive limit construction of Section \ref{sec:first-induct-limit} becomes the construction of product systems from subproduct systems as discussed in~\cite{BhMu10} and (implicitly) in~\cite{BhSk00} for Hilbert modules. A detailed discussion of inductive limits in the context of subproduct systems with manifold applications can be found in \cite{ShSk20p}.

Note that defining subproduct systems as monoidal systems in $(\Hilb^{\rm coisom},\otimes)$ gives an equivalent definition. More precisely, $\big((H_s)_{s\in\mathbb S},(\Delta_{s,t})_{s,t\in\mathbb S},\delta\big)$ is a comonoidal system in $(\Hilb^{\rm isom},\otimes)$ if and only if $\big((H_s)_{s\in\mathbb S},(\Delta_{s,t}^*)_{s,t\in\mathbb S},\delta^*\big)$ is a monoidal system in $(\Hilb^{\rm coisom},\otimes)$.

\begin{Remark}\label{sec:hilbert-spaces-1:rem:forget-HS}
 The forgetful functor $\mathcal F\colon (\FinHilb^{\rm coisom},\otimes)\to(\FinVect^{\rm surj})$ (the categories of finite-dimensional Hilbert and vector spaces, respectively) is easily seen to be a tensor functor. Tensor functors map monoidal systems to monoidal systems (just as cotensor functors do with comonoidal systems). So it follows from the previous paragraph that $(\mathcal F(H_n))_{n\in\mathbb N_0}$ yields a~standard graded algebra if $(H_n)_{n\in\mathbb N_0}$ is a~subproduct system. Similar considerations play an important role in \cite{GSk20b} where dimension sequences of finite-dimensional subproduct systems are investigated.
\end{Remark}

{\bf Additive deformations of bialgebras.}
An \emph{additive deformation} of a bialgebra $B$ is a~family of multiplication maps $\mu_t\colon B\otimes B\to B$ ($t\in\mathbb R_+$) such that $\mu_0$ is the original bialgebra multiplication, $\delta\circ\mu_t$ is pointwise continuous, $B_t:=(B,\mu_t,\Eins)$ is a unital algebra for every $t\geq0$ and $(\mu_s\otimes \mu_t)\circ (\id\otimes\tau\otimes\id)\circ(\Delta\otimes\Delta)=\Delta\circ\mu_{s+t}$ for all $s,t\geq0$, cf.\ \cite{G11} or~\cite{Wir02}. Given an additive deformation, the algebras $B_t$ with comultiplications $\Delta_{s,t}:=\Delta$ for all $s,t$ form a comonoidal system in $(\Alg_\Eins,\otimes)$, the tensor category of unital algebras with unital algebra homomorphisms and the usual tensor product. Additive deformations of $*$-bialgebras \cite{G11, Wir02} and braided ($*$)-bialgebras~\cite{GKL12} can also be defined and give rise comonoidal systems in corresponding tensor categories.

\subsection{Probability spaces}\label{sec:prob-quant-prob}

Denote by $\Prob$ the category with probability spaces as objects and measurable maps $f\colon \Omega_1\to \Omega_2$ with $\mathbb P_1\circ f^{-1}=\mathbb P_2$ as morphisms from $(\Omega_1,\mathcal{F}_1,\mathbb{P}_1)$ to $(\Omega_2,\mathcal{F}_2,\mathbb{P}_2)$.
This is a tensor category with $(\Omega_1,\mathcal{F}_1,\mathbb{P}_1)\otimes(\Omega_2,\mathcal{F}_2,\mathbb{P}_2):=(\Omega_1\times\Omega_2,\mathcal{F}_1\otimes \mathcal{F}_2,\mathbb P_1\otimes\mathbb{P}_2)$, where $\mathcal{F}_1\otimes \mathcal{F}_2$ and $\mathbb{P}_1\otimes\mathbb P_2$ are the product $\sigma$-algebra and product measure respectively. The unit object is a one-point probability space $\Lambda$. Now $\Lambda$ is clearly terminal, not initial. To use our definitions of independence and L\'evy processes we could either restate everything for the situation of the unit object being terminal, in which case we would have to reverse all arrows, or we can simply switch to the opposite category $\Prob^{\rm op}$. Both ways give basically the same definition of independence and, as noted by Franz \cite{Fra06}, this categorial independence coincides with the usual notion of stochastic independence in the following way: Let $X\colon \Omega\to E, Y\colon \Omega\to F$ be two random variables. Pushing forward the probability measure on $\Omega$ to $E$ and $F$ with $X^{-1}$ turns these into probability spaces themselves and $X$ and $Y$ can be interpreted as morphisms in $\Prob$ (or $\Prob^{\rm op}$). Now $X$ and $Y$ are stochastically independent as random variables if and only if they are categorially independent as morphisms; indeed, there $(X,Y)\colon \Omega \to E\times F$ is the unique map which makes the independence diagram commute, and it is a morphism if and only if the joint distribution $\mathbb P\circ (X,Y)^{-1}$ coincides with the product distribution $\big(\mathbb P\circ X^{-1}\big)\otimes \big(\mathbb P\circ Y^{-1}\big)$, a characterisation of stochastic independence.

Now let $(\mu_t)_{t\geq0}$ be a convolution semigroup of probability measures on the real line. In $\Prob^{\rm op}$ the addition is a~morphism from $\mathbb R$ to $\mathbb R\times\mathbb R$ and the probability spaces $A_t:=(\mathbb R,\mathcal{B},\mu_t)$ form a~comonoidal system. Indeed, $\mu_s\star\mu_t=\mu_{s+t}$ can be expressed by saying that addition transports the product measure $\mu_s\otimes \mu_t$ on $\mathbb R\times \mathbb R$ to the measure $\mu_{s+t}$ on $\mathbb R$, so addition can be interpreted as a morphism from $A_{s+r}$ to $A_s\otimes A_t$ in $\Prob^{\rm op}$. The inductive limits discussed in \Cref{theo:can-iso} are the same as the projective limit from the Daniell-Kolomogoroff theorem (taken in $\Prob$) and therefore exist.

In general, existence of projective limits in $\Prob$ is a delicate problem, see \cite{Met63} and references therein. However, it is a simple observation that the tensor product preserves projective limits. Indeed, let $\bigl((\Omega,\Sigma,\mathbb P);p_\alpha\bigr)$ and $\bigl((\Omega',\Sigma',\mathbb P');p'_\beta\bigr)$ be projective limits of projective systems $\bigl((\Omega_\alpha,\Sigma_\alpha,\mathbb P_\alpha);p_{\alpha_1,\alpha_2}\bigr)$ and $\bigl((\Omega'_\beta,\Sigma'_\beta,\mathbb P'_\beta);p_{\beta_1,\beta_2}\bigr)$, respectively. For any compatible family of morphisms $(q_\alpha,q'_\beta)\colon (X,\mathcal F,\mu)\to (\Omega_\alpha\times \Omega'_\beta,\Sigma_\alpha\otimes\Sigma'_\beta ,\mathbb P_\alpha\otimes \mathbb P'_\beta)$, there are unique $f\colon (X,\mathcal F,\mu)\to (\Omega,\Sigma,\mathbb P)$ and $f'\colon (X,\mathcal F,\mu)\to (\Omega',\Sigma',\mathbb P')$ with $p_\alpha\circ f q_\alpha$ and $p'_\beta\circ f'=q'_\beta$ by the defining property of a projective limit. The only map $F$ that makes the diagram
\[
 \begin{tikzcd}
 (\Omega\times \Omega',\Sigma\otimes\Sigma' ,\mathbb P\otimes \mathbb P')\arrow[swap]{dr}{p_\alpha\times p'_\beta}&&(X,\mathcal F,\mu)\arrow[swap]{ll}{F}\arrow{dl}{(q_\alpha,q'_\beta)}\\
 &(\Omega_\alpha\times \Omega'_\beta,\Sigma_\alpha\otimes\Sigma'_\beta ,\mathbb P_\alpha\otimes \mathbb P'_\beta)&
\end{tikzcd}
\]
commute is $F=(f,f')$, and it is readily verified that this is a morphism in $\Prob$: $\Sigma\otimes\Sigma'$ is generated by sets of the form $p_\alpha^{-1}(A)\times {p'_\beta}^{-1}(B)$ with $A\in \Sigma_\alpha,B\in\Sigma'_\beta$, and for such sets we find
\begin{gather*}
 \mu\circ(f,f')^{-1}\big(p_\alpha^{-1}(A)\times {p'_\beta}^{-1}(B)\big)\\
\qquad{} =\mu\circ (q_\alpha,q'_\beta)^{-1}(A\times B)
 =(\mathbb P_\alpha\otimes \mathbb P'_\beta)(A\times B)
 =\mathbb P_\alpha(A) \mathbb P'_\beta(B)
 =(\mathbb P\otimes \mathbb P')(A\times B).
\end{gather*}
This proves that $\bigl((\Omega\times \Omega',\Sigma\otimes\Sigma' ,\mathbb P\otimes \mathbb P');p_\alpha\times p'_\beta\bigr)$ is a projective limit of the projective system $\bigl((\Omega_\alpha\times \Omega'_\beta,\Sigma_\alpha\otimes\Sigma'_\beta ,\mathbb P_\alpha\otimes \mathbb P'_\beta);p_{\alpha_1,\alpha_2}\times p'_{\beta_1,\beta_2}\bigr)$.

\subsection{Tensor product of algebraic quantum probability spaces}

In the widest sense, an algebraic quantum probability space is a pair $(A,\varphi)$ which consists of an algebra $A$ (generalizing the algebra of bounded measurable functions on a probability space) and linear functional $\varphi\colon A\to \mathbb C$ (generalizing the expectation/integral with respect to a~probability measure). The category $\algQ$ is formed by all algebraic quantum probability spaces as objects and algebra homomorphisms $j\colon A_1\to A_2$ with $\varphi_2\circ j=\varphi_1$ as morphism from $(A_1,\varphi_1)$ to $(A_2,\varphi_2)$. By $\algQ_\Eins$ we denote the subcategory whose objects consist of a unital algebra with a unital linear functional and whose morphisms are only those morphisms of $\algQ$ which are unital algebra homomorphisms.

Maybe the simplest tensor category in this context is $(\algQ_\Eins,\otimes)$, with tensor product $(A_1,\varphi_1)\allowbreak \otimes(A_2,\varphi_2):=(A_1\otimes A_2,\varphi_1\otimes \varphi_2)$ and unit object $(\mathbb C,\id_\mathbb C)$. The unit object is initial in $\algQ_\Eins$, so we have a derived notion of independence for morphisms. Two morphisms $j_i\colon(A_i,\varphi_i)\to(A,\Phi)$, $i\in\{1,2\}$, are independent if and only if the images of $j_1$ and $j_2$ commute with each other and $\Phi\circ\mu_{A}\circ(j_1\otimes j_2)=(\Phi\circ j_1)\otimes(\Phi\circ j_2)$. Indeed, if there exists an independence morphism $h\colon (A_1,\varphi_1)\otimes\,(A_2,\varphi_2)\to(A,\Phi)$, then it must hold that
\[
h(a\otimes b)=h(a\otimes\Eins)h(\Eins\otimes b)=j_1(a)j_2(b)=\mu_A \circ(j_1\otimes j_2)(a\otimes b)
\]
and
\[
\Phi\circ\mu_A \circ(j_1\otimes j_2)=\Phi\circ h=\varphi_1\otimes\varphi_2=(\Phi\circ j_1)\otimes\,(\Phi\circ j_2).
\]
The fact that $h$ is an algebra homomorphism implies
\[
j_1(a)j_2(b)=\mu_A\circ(j_1\otimes j_2)(\Eins\otimes b\cdot a\otimes \Eins)=j_2(b)j_1(a)
\]
so the images commute. On the other hand, if $\Phi\circ\mu_{A}\circ(j_1\otimes j_2)=(\Phi\circ j_1)\otimes(\Phi\circ j_2)$, then $h:=\mu_A\circ(j_{1}\otimes j_2)$ respects the linear functionals, and if the images commute, it is an algebra homomorphism, hence a (unique) independence morphism.

Now let $B$ be a bialgebra with comultiplication $\Delta$ and counit $\delta$. Then linear functionals on $B$ can be convolved. Every convolution semigroup $(\varphi_t)_{t\geq0}$ now gives rise to a comonoidal system $A_t=(B,\varphi_t)$, because $\varphi_{s+t}=\varphi_s\star\varphi_t$ is equivalent to $\Delta$ being a morphism in $\algQ_\Eins$ from $A_{s+t}$ to $A_s\otimes A_t$ and $\varphi_0=\delta$ is equivalent to $\delta$ being a morphism from $A_0$ to $(\mathbb C,\id_\mathbb C)$. If one also requires pointwise continuity, the resulting L\'evy processes are the quantum L\'evy processes on bialgebras whose theory is developed in \cite{Sch93}.

The setting of comonoidal system gives us more freedom. It is not necessary that the algebras of the $A_t$ are all the same. This situation comes up when one considers L\'evy processes on additive deformations of bialgebras as discussed above. Again, a convolution semigroup of linear functionals $\varphi_t$ on the algebras $B_t:=(B,\mu_t)$ gives rise to a comonoidal system $A_t=(B_t,\varphi_t)$ and corresponding L\'evy processes.

\subsection{Universal products}\label{sec:universal-products}

Independence in quantum probability is usually implemented by a \emph{universal product}, which is a~prescription $\boxdot$ that assigns to two linear functionals on algebras $\mathcal A_1$, $\mathcal A_2$ a new linear functional $\varphi_1\boxdot\varphi_2$ on the \emph{free product} $\mathcal A_1\sqcup\mathcal A_2$ such that the bifunctor $((\mathcal A_1,\Phi_1),(\mathcal A_2,\Phi_2))\mapsto(\mathcal A_1\sqcup\mathcal A_2$, $\Phi_1\boxdot \Phi_2)$ turns the category $\algQ$ of quantum probability spaces into a tensor category with the canonical embeddings $\mathcal A_i\hookrightarrow\mathcal A_1\sqcup\mathcal A_2$ as inclusions, see, e.g., \cite{bGSc02}. Categorial independence of $j_i\colon (\mathcal B_i,\psi_i)\to (\mathcal A,\Phi)$ in this situation means that the \emph{noncommutative joint distribution} $\Phi\circ (j_1\sqcup\dots\sqcup j_n)$ coincides with the universal product of marginal distributions $\psi_1\boxdot\dots\boxdot\psi_n=(\Phi\circ j_1)\boxdot\dots \boxdot (\Phi\circ j_n)$.

Note that this is somehow only a slight variation of the principle which led to independences in the algebraic situations. We again have an underlying category with coproduct, the category of algebras arbitrary algebra homomorphisms, the coproduct being the free product. Instead of going directly to a subcategory, we enrich the objects by additional structure (the linear functionals), and restrict to morphisms which respect the additional structure.

In quantum probability there are five well known notions of independence which come from universal products, namely tensor independence, freeness, boolean independence, monotone, and antimonotone independence \cite[Section 1.8]{FrSk16}.

 \begin{Remark}
 In order to motivate the definition of universal products a bit more, suppose that there is any bifunctor $\boxtimes$ which turns $(\algQ,\boxtimes)$ into a tensor category such that the unit object is the initial object $\{0\}$. By \Cref{theo:compincl-uoinitial}(b), there are inclusions $Q_i\to Q_1\boxtimes Q_2$. Denote by $A_Q$ and $\varphi_Q$ the algebra and the linear functional of $Q$, i.e., $Q=(A_Q,\varphi_Q)$. Then the inclusions are algebra homomorphisms $A_{Q_i}\to A_{Q_1\boxtimes Q_2}$. Denote by $A_1\sqcup A_2$ the free product of the algebras~$A_1$ and~$A_2$. The universal property of the free product establishes an algebra homomorphism $\iota^{1}\sqcup\iota^{2}\colon A_{Q_1}\sqcup A_{Q_2}\to A_{Q_1\boxtimes Q_2}$. For linear functionals $\varphi_i\colon A_i\to \mathbb C$ put $Q_i:=(A_i,\varphi_i)$ and define $\varphi_1\odot \varphi_2\colon A_1\sqcup A_2$ as $\varphi_1\odot \varphi_2:=\varphi_{Q_1\boxtimes Q_2}\circ\big(\iota^{1}\sqcup\iota^{2}\big)$. Then $\odot$ fulfills the axioms of a universal product (cf.\ \cite[Definition~3.1]{GLa15}).
 \end{Remark}

A dual semigroup is by definition a comonoid in the category $\Alg$ of algebras. With respect to a fixed universal product ${\boxdot}$, a convolution product for linear functionals on a dual semigroup can be defined via $\varphi_1\star\varphi_2:=(\varphi_1\boxdot\varphi_2)\circ\Delta$, where $\Delta\colon D\to D\sqcup D$ is the comultiplication of the dual semigroup $D$. As in the tensor product case, convolution semigroups can be studied and lead to comonoidal systems $A_t:=(D,\varphi_t)$ and finally to L\'evy processes whose increments are independent in the sense given by the universal product (for example free, monotone or Boolean independence).

Universal products in the category $\algQ$ have been almost completely classified (see \cite{bGSc02,GLa15,Mur03, Spe97}). Correspondingly, the theory of L\'evy processes can be (and has been) dealt with by studying the special cases~\cite{bGSc05}. There are two important generalizations of universal products for which a classification is out of reach. First, Bozeiko, Leinert and Speicher studied so-called \emph{c-freeness}, which is obtained by taking products of pairs of linear functionals on an algebra~\cite{BLS96}. Following this, Hasebe introduced the \emph{indented product}, which is a~product for triples of linear functionals~\cite{Has10p} and studied also a product for infinitely many functionals~\cite{Has11}. So instead of~$\algQ$, these products give bifunctors for a category $\algQ_n$ whose objects are tuples $(A,\varphi_1,\dots,\varphi_n)$. In another direction, in a series of papers~\cite{Voi14,Voi16a,Voi16b} Voiculescu presented a~fascinating new notion of independence which he calls \emph{bifreeness}. To understand bifreeness as categorial independence, algebras have to come with a free product decomposition $A=A_l\sqcup A_r$ into two \emph{faces}, where $A_l$ contains the \emph{left variables} and $A_r$ the \emph{right variables}.\footnote{In practice, left and right variables are not always in free product position from the beginning. However, for the purpose of describing the independence, one can always artificially form the free product of the algebras of left and right variables.} One comes to study the category $\algQ^{m}$ whose objects are tuples $(A, A_1,\dots, A_m,\varphi)$, $\varphi\colon A=A_1\sqcup\cdots\sqcup A_m\to \mathbb C$. It is even possible to combine the two generalizations and work in a category~$\algQ_n^{m}$, defined in the obvious way. For all of these categories one can define (multivariate) universal products, and develop the theory of L\'evy processes, as discussed in this paper, as well as other topics related to independence, as has been successfully done with cumulants in~\cite{MaSc17}. In the last years, such multivariate independences have received a lot of attention and many more examples have been found and studied to some extent \cite{G17p,G21p,GHU21p,GHS20,GuSk17,GuSk19,Liu18p,Liu19}. The richness of examples when compared to the univariate case underlines the value of general methods applying to all multivariate universal independences at once.

\section{Outlook}

In this article we defined independence and L\'evy processes in a synthetic probability setting and we could establish a general Daniell--Kolmogorov type reconstruction of categorial L\'evy processes from its comonoidal system of `marginal distributions'. We assume that the following directions of research could be very fruitful continuations of this work.
\begin{itemize}\itemsep=0pt
\item In classical and quantum probability, L\'evy processes are often seen as the basic buiding blocks from which more complicated (in particular nonstationary) processes are composed by means of (quantum) stochastic calculus. Can categorial L\'evy processes in the same way be used to construct a larger, but still controllable class of processes in a synthetic probability setting?
\item In many cases where one deals with L\'evy processes (classical stochastic processes on topological groups, quantum stochastic processes on bialgebras or dual groups) some weak continuity is assumed, so that the marginal distributions forn a differentiable semigroup which is determined by a generator. Is it possible to implement such ideas in the framework of categorial L\'evy processes?
\item Are there general contexts in which it is possible to obtain a canonical decomposition of categorial L\'evy processes, like the L\'evy Khintchine decomposition into a~Gaussian and a~Poisson part, or even a~classification comparable to the L\'evy Khintchine formula or Hunt's formula?
\end{itemize}

\subsection*{Acknowledgements}
The authors are grateful to Uwe Franz, Michael Skeide, Tobias Fritz, Simeon Reich, Orr Shalit and the anonymous referees for helpful comments on earlier drafts of this article.
The work of MG and MS was supported by the German Research Foundation (DFG), project number 397960675. MG's work was carried out partially during the tenure of an ERCIM `Alain Bensoussan' Fellowship Programme.

%\bibliographystyle{sigma}
%\bibliography{sigma22-07Gerhold}

\pdfbookmark[1]{References}{ref}
\LastPageEnding

\end{document}